\documentclass[12pt]{article}
\usepackage{amsmath, amssymb}
\usepackage{amsfonts}
\usepackage[dvips]{graphicx}

\newtheorem{theorem}{Theorem}[section]
\newtheorem{lemma}{Lemma}[section]
\newtheorem{definition}{Definition}[section]

\def\QED{\mbox{\rule[0pt]{1.5ex}{1.5ex}}}

\def\endproof{\hspace*{\fill}~\QED\par\endtrivlist\unskip}
\def\keywords{\vspace{-.3em}
    \if@twocolumn
      \small\it Keywords\/\bf---$\!$%
    \else
      \begin{center}\small\bf Keywords\end{center}\quotation\small
    \fi}
\def\endkeywords{\vspace{0.6em}\par\if@twocolumn\else\endquotation\fi
    \normalsize\rm}
\def\Label#1{\label{#1}\ [\ #1 \ ]\ }
\def\Label{\label}

\title{Differentiability of eigenfunctions of the closures \\ of differential operators \\ with rational coefficient functions}
\author{Fuminori SAKAGUCHI\thanks{Faculty of Engineering, University of Fukui, 3-9-1 Bunkyo, Fukui 910-8507, Japan ({\tt fsaka@u-fukui.ac.jp})} \and Masahito HAYASHI\thanks{Graduate School of Information Sciences, Tohoku University, Sendai 980-8579, Japan ({\tt hayashi@math.is.tohoku.ac.jp}) Centre for Quantum Technologies, National University of Singapore, 3 Science Drive 2, Singapore 117542}}

\date{}

\begin{document}
\maketitle

\begin{abstract}
In this paper, for an operator defined by the action of an $M$-th order differential operator with rational-type coefficients on the function space 
$L_{(k_0)}^2(\mathbb{R} ):=\{ f:{\rm measurable} | \,\, \|f\|_{(k_0)} <\infty  \}$ 
with norm  $\|f\|_{(k_0)}^2:= \int  |f(x)|^2  (x^2+1)^{k_0}  dx$\,   ($k_0\in\mathbb{Z} $), 
we prove the regularity (continuity and differentiability up to $M$ times) of the eigenfunctions of its closure (with respect to the graph norm), 
except at singular points of the corresponding ordinary differential equation
without any assumptions for the Sobolev space, i.e., 
without any assumptions about the $m$-th order derivatives of the eigenfunctions with $m=1,2,\ldots M-1$. 
(For the special case of $k_0=0$, we prove this regularity for the usual $L^2(\mathbb{R} )$.)
Especially, we show a one-to-one correspondence between the eigenfunctions of its closure and the solutions in $C^M(\mathbb{R} )\cap L_{(k_0)}^2(\mathbb{R} )$ 
of the corresponding differential equation under the condition above 
when there is no singular point for this differential equation.
This one-to-one correspondence is shown in the basic framework of 
an algorithm proposed in our preceding paper, 
which can determine all solutions in $C^M\cap L_{(k_0)}^2(\mathbb{R} )$ 
of the ordinary differential equation then.
% In this framework,
% the differential operator is treated as an operator from 
% a Hilbert space ${\cal H}$ to another Hilbert space ${\cal H}^\Diamond $,
% and it can be represented in matrix form with  appropriate basis systems of 
% ${\cal H}$ and ${\cal H}^\Diamond $.
% In accordance with this matrix representation, we transform eigenfunctions in ${\cal H}$ 
% to square-summable number sequences satisfying the matrix-vector equation.
% The truncation of this square-summable number sequence
% yields an appropriate approximation of the eigenfunction by an $M$-times differentiable function.
% We can show that the eigenfunction belongs to $C^M(\mathbb{R} )$
% when the pair of Hilbert spaces ${\cal H}$ and ${\cal H}^\Diamond $
% satisfies several conditions and this approximation has point-wise convergence.
\end{abstract}

AMS: 35B65, 34L99, 52C05

\section{Introduction}
\Label{sec:in}
When we treat the eigenfunction problem of the closure of an $M$-th order differential operator on a Hilbert space with a certain boundary condition,
we should be careful to distinguish this problem from the problem of finding solutions in the space of $M$-times differentiable functions $C^M(\mathbb{R} )$ of the differential equation 
described by this differential operator definable only in $C^M(\mathbb{R} )$. If a solution to the latter problem belongs to the Hilbert space above and satisfies the boundary condition, it is always an eigenfunction of the former problem from the definition. However, it is not necessarily the case that eigenfunctions of the former problem belong to $C^M(\mathbb{R} )$. Hence, the regularity (continuity and differentiability up to $M$ times) of the eigenfunctions of the former problem should be examined carefully.

In the theory of elliptic operators~\cite{Gil}, this problem has been discussed under assumptions %about 
for  
the Sobolev space, i.e., the assumption that the $m$-th order derivatives of the eigenfunction with $m=1,2,\ldots M-1$ belong to $L^p$-space. These assumptions are often required for the validity of numerical methods that solve differential equations by projection to finite dimensional  subspaces (Ritz-Galerkin and Petrov-Galerkin 
methods~\cite{KrVa}~\cite{Bru}, for example).

On the other hand, in this paper, for a class of Hilbert spaces containing $L^2(\mathbb{R} )$, we will discuss the regularity problem above under 
%a certain condition, 
several conditions,  
without any assumptions concerning  the $m$-th order derivatives of the eigenfunction. 
The condition in our discussion is that the differential operator has rational coefficient functions 
and its  
characteristic equation (eigenvalue equation) has no singular point. Under this condition, we prove that the regularity above is always guaranteed. 

This discussion  can be generalized, even for a differential operator (with rational coefficient functions) whose  
characteristic equation has singular points, by excluding only the regularity at the singular points of the characteristic equation. Especially for Fuchsian-type differential operators, we give a stronger statement than general cases.
% The condition in our discussion is that the differential operator has polynomial-type coefficient functions such that the coefficient function of the highest order has no zero point. This condition can be generalized, even for a differential operator with rational coefficient functions such that the denominators have no zero points and the coefficient function of the highest order also has no zero point. Under this condition, we prove that the regularity above is always guaranteed.

% This discussion can be generalized for the cases where the above condition for the highest-order coefficient function and the denominators is not necessarily satisfied, by excluding only the regularity at the zero points of the highest-order coefficient function and the denominators.

The proof is based on a one-to-one correspondence between the `regular' solutions in the Hilbert space of the differential equation and the square-summable number-sequence solutions of simultaneous linear equations described by a kind of matrix representation of the action of the differential operator, which is guaranteed under several conditions. In this paper, we will clarify how we can show regularity using this one-to-one correspondence. 

The contents of this paper are as follows: 
Section \ref{sec:pp} introduces the basic framework used for the proof. 
Firstly, in Subsection \ref{subsec:s21} we clarify precisely what has to be proved. 
Next, Subsection \ref{subsec:s22} provides a more general framework in which the regularity problem can be discussed, and it shows the conditions that are required 
%within 
for the base of  
this framework. 
In Subsection \ref{subsec:s24}, 
In order to treat the argument given in Subsection \ref{subsec:s21},
we give a concrete structure for the general theory given in Subsection \ref{subsec:s22}.
Subsection \ref{subsec:s25} explain how to 
apply general theory given in Subsection \ref{subsec:s22} to the concrete structure given in Subsection \ref{subsec:s24} for showing the statement given in Subsection \ref{subsec:s21}.
Section \ref{sec:sm} is devoted to a proof of a property mentioned in 
Subsection \ref{subsec:s24}.
Section \ref{sec:nenps} is devoted to a proof of a theorem mentioned in 
Subsection \ref{subsec:s22}.

\section{Basic framework of this paper}
\Label{sec:pp}

\subsection{`Regularity' of eigenfunctions to be shown}
\Label{subsec:s21} 
In this subsection, we rigorously describe 
the regularity problem to be solved in this paper.
In this paper, we treat the differential operator 
\begin{eqnarray}
{R}(x,{\textstyle \frac{d}{dx}})
:=\sum_{m=0}^M {r}_m (x) \bigl({\textstyle \frac{d}{dx}}\bigr)^m  \, 
\Label{eqn:def-diff-op}
\end{eqnarray}
on the space of $M$-times differentiable functions $C^M(\mathbb{R} )$.

In order to treat 
the ODE $ {R}(x,\frac{d}{dx})f(x)=\lambda f(x) $ 
using functional analysis,
we have to define the differential operator in a complete function space.

In the present paper, 
we focus on the function space
$L_{(k_0)}^2(\mathbb{R} )$, which is defined by 
\begin{eqnarray}
L_{(k_0)}^2(\mathbb{R} ):=\{ f:{\rm measurable}\, | \, \|f\|_{(k_0)} <\infty  \}
\Label{eqn:Ha1}\end{eqnarray}
with inner product $\displaystyle \langle f,g\rangle _{(k_0)}= \int_{-\infty }^\infty f(x) \overline{g(x)} (x^2+1)^{k_0}  dx$ and norm  
$\displaystyle \|f\|_{(k_0)}^2
= \int_{-\infty }^\infty |f(x)|^2 (x^2+1)^{k_0}  dx$,
parametrized by $k_0\in\mathbb{Z} $. 
Here, in the special case of $k_0=0$,\,  $L_{(0)}^2(\mathbb{R} )$ is identical to the usual $L^2(\mathbb{R} )$. 
Then, 
the operator $\widetilde{A}_{R,L_{(k_0)}^2(\mathbb{R} )}$ is defined by the action of ${R}(x,{\textstyle \frac{d}{dx}})$ with domain
\begin{eqnarray}
\hspace{5mm}
D(\widetilde{A}_{R,L_{(k_0)}^2(\mathbb{R} )}):=
\{f\in C^M(\mathbb{R} ) 
\cap L_{(k_0)}^2(\mathbb{R} ) \, |\, {R}(x,{\textstyle \frac{d}{dx}}) f \in L_{(k_0)}^2(\mathbb{R} ) \},\,   
\Label{eqn:def-op-A}\end{eqnarray}
and its closure $A_{R,L_{(k_0)}^2(\mathbb{R} )}$ with respect to the graph norm~\cite{ReSi1}. 

In general, 
an eigenfunction of the closed extension of the given differential operator 
does not necessarily yield a solution of 
the ODE $ {R}(x,\frac{d}{dx})f(x)=\lambda f(x) $.
This is because there is a possibility that 
the eigenfunction is not an $M$-times differentiable function.
This problem is called %a regularity problem for a differential operator.
the regularity problem for a differential operator.

In the present paper, we prove that the eigenfunction of the operator 
$A_{R,L_{(k_0)}^2(\mathbb{R} )}$ always does yield a
solution of the ODE $ {R}(x,\frac{d}{dx})f(x)=\lambda f(x) $.
Since this problem depends on the singularity of the differential equation,
we need the following definition.
A real number $x \in \mathbb{R}$ is called a singular point of
a differential operator 
${R}(x,{\textstyle \frac{d}{dx}})$
if $x$ is a singular point of the differential equation 
${R}(x,{\textstyle \frac{d}{dx}})f(x)=\lambda f(x)$ for some real number $\lambda$.
Indeed, the above definition does not depend on $\lambda$.
When the coefficient functions of 
${P}(x,{\textstyle \frac{d}{dx}}):=\sum_{m=0}^M {p}_m (x) \bigl({\textstyle \frac{d}{dx}}\bigr)^m
$ are polynomial,
the set of its singular points equals the set of zero points of
its coefficient function ${p}_M$ of the highest degree, which is written as 
${p}_M^{-1}(0)$.
When the coefficient functions 
${r}_m(x)$ $(m=0,1,\ldots M)$
of ${R}(x,{\textstyle \frac{d}{dx}})$ are rational functions,
we denote 
the least common multiple of the denominators of ${r}_m(x)$
by $l(x)$.
Then,
the set of its singular points equals the set of zero points of
$l\cdot {r}_M(x):=l(x) {r}_M(x)$, which is written as 
$(l\cdot{r}_M)^{-1}(0)$.

First, we obtain the following theorem.

\begin{theorem}
\Label{thm:main}
Assume that 
a differential operator 
${R}(x,{\textstyle \frac{d}{dx}})$ has no singular points 
and 
its coefficient functions ${r}_m(x)$ $(m=0,1,\ldots M)$ are rational functions.
Then,
an element  $f\in L_{(k_0)}^2(\mathbb{R} )$
is an eigenfunction of $A_{R, L_{(k_0)}^2(\mathbb{R} )}$ with the eigenvalue $\lambda$
if and only if
$f$ belongs to $C^M(\mathbb{R}) \cap L_{(k_0)}^2(\mathbb{R} )$
and satisfies the ODE ${R}(x,{\textstyle \frac{d}{dx}}) f(x)=\lambda f(x)$.
\end{theorem}
This theorem is one of the main statements to be proved in this paper. 
When the differential operator 
$R(x,{\textstyle \frac{d}{dx}})$ has singular points,
the following extension holds.

\begin{theorem}
\Label{thm:main2}
When 
the coefficient functions ${r}_m(x)$ $(m=0,1,\ldots M)$ 
of a differential operator 
${R}(x,{\textstyle \frac{d}{dx}})$ 
are rational functions,
then
any eigenfunction of $A_{R, L_{(k_0)}^2(\mathbb{R} )}$ belongs to $C^M(\mathbb{R} \setminus (l\cdot{r}_M)^{-1}(0))$
for any integer $k_0$. 
\end{theorem}

In the next subsection, we will give a more general argument, which includes 
Theorem \ref{thm:main} as a special case.

\subsection{Regularity in a more general framework}
\Label{subsec:s22}
In this subsection, 
we treat the regularity problem  
%in a general Hilbert function space ${\cal H}$ on a real line $\mathbb{R}$.
in a general Hilbert space ${\cal H}$ of functions on the real line $\mathbb{R}$.
That is, we give three conditions equivalent to the solution of 
the ODE $ {P}(x,\frac{d}{dx})f(x)=\lambda f(x)$ 
in a general Hilbert function space ${\cal H}$, where %i.e., 
we convert the ODE to square-summable solutions of a matrix-vector equation (simultaneous linear equations) defined in the following general framework.

Now, we introduce another general Hilbert function space ${{\cal H}^\Diamond}$
as a Hilbert function space on $\mathbb{R}$ which contains (as a subset) the original Hilbert function space ${\cal H}$. 
In general, the inner product of 
${\cal H}$ 
is distinct from the inner product of 
${\cal H}^\Diamond$, 
%and hence the space ${\cal H}^\Diamond$ is not a subspace of  ${\cal H}$ 
whereas ${\cal H}$ is a subset of ${\cal H}^\Diamond$. 
By treating the differential operator as an operator from ${\cal H}$ to ${{\cal H}^\Diamond}$,
we are able to utilize a `matrix representation' of the ODE with respect to  appropriate basis systems.
The key point of the method that we present is  
the difference between the inner products of %both spaces ${\cal H}$ and ${{\cal H}^\Diamond}$.
the spaces ${\cal H}$ and ${{\cal H}^\Diamond}$.

Define the operator $\widetilde{A}_{P,{\cal H}}$ as the action of ${P}(x,{\textstyle \frac{d}{dx}})$ with domain
\begin{eqnarray}
D(\widetilde{A}_{P,{\cal H}}):= \{ f\in 
 C^M(\mathbb{R} ) 
\cap {\cal H} \, |\, {P}(x,{\textstyle \frac{d}{dx}}) f \in {\cal H} \},\,   
\Label{eqn:def-op-A_2}\end{eqnarray}
and its closure $A_{P,{\cal H}}$ with respect to the graph norm. 
Next, we introduce an operator from ${\cal H}$ to ${{\cal H}^\Diamond}$.
Define the operator $\widetilde{B}_{P,\lambda ,   {\cal H}, {{\cal H}^\Diamond}}$ as the action of 
\begin{eqnarray}
P(\lambda; x,{\textstyle \frac{d}{dx}}):={P}(x,{\textstyle \frac{d}{dx}})-\lambda 
I  
\,\,\,\,\,\,\,\,\,\,\,\,\,(I:\mbox{identity op.})
\Label{eqn:def-diff-op_2}\end{eqnarray}
with domain
\begin{eqnarray}
D(\widetilde{B}_{P,\lambda ,   {\cal H}, {{\cal H}^\Diamond}}):=
\{ f\in 
C^M(\mathbb{R} ) 
\cap {\cal H} \, |\, P(\lambda; x,{\textstyle \frac{d}{dx}}) f \in {{\cal H}^\Diamond} \},\,   
\Label{eqn:def-op-p-HH}\end{eqnarray}
and its closure $B_{P,\lambda , {\cal H}, {{\cal H}^\Diamond}}$ with respect to the corresponding graph norm $\|\cdot\|_{\cal H}+ \|P(\lambda; x,\frac{d}{dx})\,\cdot \|_{{{\cal H}^\Diamond}}$.

%The main result is the equivalence between 
%the solutions in the `matrix representation' 
% the solutions of the system of simultaneous linear equations corresponding to the matrix representaion $b_{m}^n:=\langle B_{P,\lambda , {\cal H}, {{\cal H}^\Diamond}} e_n,\, e^\Diamond_m \rangle_{{{\cal H}^\Diamond}}$  
% and the solutions of 
% the ODE $ {P}(x,\frac{d}{dx})f(x)=\lambda f(x)$ in a general Hilbert function space ${\cal H}$
% under conditions {\bf C1-C3}, and {\bf C2.1-C2.4}  below:

In order to using a band-diagonal structure in 
the close operator $B_{P,\lambda , {\cal H}, {{\cal H}^\Diamond}}$,
we introduce
Conditions {\bf C1-C3}, {\bf C1$^{+}$}, {\bf C2$^{+}$}, and {\bf C2.1-C2.3} 
for the quintuplet consisting of
the linear differential operator $P(\lambda; x,{\textstyle \frac{d}{dx}})$,
the Hilbert spaces ${\cal H}$ and ${{\cal H}^\Diamond }$,
and their CONSs 
$\{e_n \, \}_{n=0}^{\infty}$ and 
$\{e_n^\Diamond \}_{n=0}^{\infty}$,
which is abbreviated to 
$(P(\lambda; x,{\textstyle \frac{d}{dx}}),{\cal H},
\{e_n \, \}_{n=0}^{\infty},
{{\cal H}^\Diamond },
\{e_n^\Diamond \}_{n=0}^{\infty})$.
In what follows, $\langle \cdot , \, \cdot \rangle _{{{\cal H}^\Diamond }}$
and 
$\langle \cdot ,
 \, \cdot \rangle _{{\cal H}}$ denote the inner products of ${{\cal H}^\Diamond }$ and ${\cal H}$ respectively.
These conditions are shown to hold in several examples for $P(x,{\textstyle \frac{d}{dx}})$ later.

\begin{description}
\item[C1]
%The CONS $\{e_n \, \}_{n=0}^{\infty}$ of ${\cal H}$ satisfies
For any $n$, $e_n$ belongs to $D(\tilde{B}_{P,\lambda,{\cal H},{{\cal H}^\Diamond}})$.
\item[C1$^+$]
There exists a positive function $\upsilon$ in $C^M(\mathbb{R} )$ such that $\displaystyle \langle f,g \rangle _{\cal H}=\int_{-\infty }^\infty \!\!\! f(x)\overline{g(x)}\upsilon(x)dx$.
\item[C2]
There exists an integer $\ell _0 $ 
such that $b_{m}^n:=\langle 
B_{P,\lambda,{\cal H},{{\cal H}^\Diamond}}
 e_n, e_m^\Diamond \rangle _{{{\cal H}^\Diamond }}=0$ 
when $|n-m|> \ell _0 $.
\item[C2$^+$]
There exists a positive function ${\upsilon^\Diamond}$ in $C^M(\mathbb{R} )$ such that  $\displaystyle \langle f,g \rangle _{{\cal H}^\Diamond}
=\int_{-\infty }^\infty \!\!\! f(x)\overline{g(x)}{\upsilon^\Diamond}(x)dx$.
% \item[C1.1]
% If a number sequence $\{f_n\}_{n=0}^\infty \in \ell ^2(\mathbb{Z}^+)$ satisfies $\, \displaystyle \lim _{N\to\infty } \Bigl\| \Bigl(\sum _{n=0}^N f_n e_n \Bigr) -f\, \Bigr\|_{\cal H}=0$ with $f\in C^M(\mathbb{R} )$, \, then $\, \displaystyle \lim _{N\to\infty }\sum _{n=0}^N f_n e_n(x) =f(x)$ holds for any $x\in\mathbb{R}$.
\item[C2.1]
$\displaystyle \sup_{n\in \mathbb{Z}^+\backslash\{ 0\} }\frac{|b_m^n|}{n^M}< \infty $.
%, where $M$ is  the order of ${P}(x,\frac{d}{dx})$.   
%\item[C2.3]
%For any fixed $a$ and $b$ s.t. $-\infty < a < b <\infty$, if a function $f$ in $C^M(\mathbb{R} )$ satisfies $f(x)=0$ for $x\le a$ and $f^{(m)}(b)=0$ for $m=0,1,\ldots ,M-2$, then $\displaystyle \sup_{n\in\mathbb{Z}^+} n^M \left|\int_a^b f(x) \,e^\Diamond_n(x) dx \right|< \infty$.
\item[C2.2]
The basis functions $e^\Diamond_n$ $(n\in \mathbb{Z}^+ )$ belong to $C^M(\mathbb{R} )$ and there exists a first-order differential operator $N(x,\frac{d}{dx})=n_1(x)\frac{d}{dx}+n_0(x)$ satisfying (a) and (b) below:

(a): The functions $n_1$ and $n_0$ belong to  $C^{M-1}(\mathbb{R} )$ 

(b): There exist real numbers $\lambda _n$ $(n\in\mathbb{Z}^+ )$ such that  $N(x,\frac{d}{dx})\,e^\Diamond_n=\lambda _ne^\Diamond_n$  for any  $ n\in\mathbb{Z}^+ $,\,   and  $\displaystyle \liminf_{n\to\infty }\frac{|\lambda _n|}{n}>0$.
\item[C2.3]
There exists a function $\tilde{a}$ in $C^0(\mathbb{R} )$ such that $^\forall n\in\mathbb{Z}^+ $ and $ ^\forall x\in \mathbb{R} $, $|e^\Diamond_n(x)|\le \tilde{a}(x)$.
\item[C3]
There exists a linear operator $C_{P,\lambda,{\cal H},{{\cal H}^\Diamond}}$ 
with domain $D(C_{P,\lambda,{\cal H},{{\cal H}^\Diamond}})$ 
from a dense subspace of ${{\cal H}^\Diamond }$ to ${\cal H}$
such that $e_m^\Diamond \in D(C_{P,\lambda,{\cal H},{{\cal H}^\Diamond}})$ and 
$\langle B_{P,\lambda,{\cal H},{{\cal H}^\Diamond}} f, e_m^\Diamond \rangle _{{{\cal H}^\Diamond }}=\langle f, 
C_{P,\lambda,{\cal H},{{\cal H}^\Diamond}} e_m^\Diamond \rangle _{{\cal H}}$
for $f \in D(\tilde{B}_{P,\lambda,{\cal H},{{\cal H}^\Diamond}})$.
\end{description}

Our main issue is the correspondence between the following two kinds of solutions 
One kind of solutions are 
the 
 square-summable  
solutions of the system of simultaneous linear equations corresponding to the matrix 
representation 
$b_{m}^n:=\langle B_{P,\lambda , {\cal H}, {{\cal H}^\Diamond}} e_n,\, e^\Diamond_m \rangle_{{{\cal H}^\Diamond}}$. The other kinds of solutions are the solutions of 
the ODE $ {P}(x,\frac{d}{dx})f(x)=\lambda f(x)$ in a general Hilbert function space ${\cal H}$. 

Due to the condition {\bf C3},
the basis $e^\Diamond_m$ belongs to the domain of the adjoint operator $B_{P,\lambda , {\cal H}, {{\cal H}^\Diamond}}^{\, *}$. In the following two conditions, $M$ denotes the order of ${P}(x,\frac{d}{dx})$.

%By means of a pair of orthonormal basis systems $\{e_n\, |\, n\in\mathbb{Z}^+ \} $ and $\{e^\Diamond_n\, |\, n\in\mathbb{Z}^+ \} $ of ${\cal H}$ and ${{\cal H}^\Diamond}$, respectively, we introduce the `matrix representation':
% \begin{eqnarray}
% b_m^n:=\langle B_{P,\lambda , {\cal H}, {{\cal H}^\Diamond}}\,  e_n,\,  e^\Diamond_m \rangle_{{{\cal H}^\Diamond}} \, , 
% \Label{eqn:def_matrix_elem}\end{eqnarray}
% for $B_{P,\lambda , {\cal H}, {{\cal H}^\Diamond}}$ defined above, where $\langle \cdot ,\cdot\rangle_{{{\cal H}^\Diamond}}$ denotes the inner product of ${{\cal H}^\Diamond}$. 
%where the parameters are omitted in the notations for simplicity. By these representations, we will introduce the following number-sequence space:

With $b_m^n$ defined in {\bf C2}, define the solution space $V$ as a space of number sequences 
\begin{eqnarray}
\displaystyle V:=\bigl\{ \{f_n\}_{n=0}^\infty  \, | \, \sum_{n=0}^\infty b_m^n f_n =0  \,\,\, (m\in\mathbb{Z}^+)\bigr\} .
\Label{eqn:sp-number-seq}\end{eqnarray}
With this definition, one of the `equivalent conditions' mentioned above is 
$\{f_n\}_{n=0}^\infty\in V\cap \ell ^2(\mathbb{Z}^+)$. 
As is shown later, the following theorem holds.

\begin{theorem}
%\Label{thm:equivalence_3spaces}
\Label{thm:equivalence1}
When 
the quintuplet
$(P(x,{\textstyle \frac{d}{dx}}),{\cal H},
\{e_n \, \}_{n=0}^{\infty},
{{\cal H}^\Diamond },
\{e_n^\Diamond \}_{n=0}^{\infty})$
satisfies 
{\bf C1-C3}, {\bf C2$^+$}, and {\bf C2.1-C2.3},
then the relations 
$
(i)\Rightarrow
(ii)\Rightarrow
(iii)\Rightarrow
(iv)$
holds concerning the conditions for $f \in {\cal H}$:

$(i)$: \, $f\in {\rm dom} \, A_{P,{ \cal H}}$ and  
$A_{P,{ \cal H}}\, f = \lambda f$.

$(ii)$: \, $f\in {\rm dom} \, 
B_{P,\lambda , {\cal H}, {{\cal H}^\Diamond}}$ and  $B_{P,\lambda , {\cal H}, {{\cal H}^\Diamond}}\, f = 0$

$(iii)$: \, $
f\in{\cal H} \mbox{ and } 
\{\langle f, e_n \rangle_{\cal H}\}_{n=0}^\infty \in V \cap \ell ^2(\mathbb{Z}^+)$.

$(iv)$: $f\in 
C^M(\mathbb{R}\setminus S ) 
\cap {\cal H}$ and 
$^\forall x \in \mathbb{R}\setminus S$, 
$ {P}(x,\frac{d}{dx})f(x)=\lambda f(x) $,
$S$ is the set of singular points of $(P(x,{\textstyle \frac{d}{dx}})$.
\end{theorem}

We can prove a stronger argument than the above theorem
with an additional condition.

In this paper, we define the following two types of Fuchsian operators: 
\begin{definition}
$P(x,\frac{d}{dx})$ is called `Fuchsian of Type I' 
if its all singular points are 
regular singular points of 
the ODE $P(x,\frac{d}{dx})f(x)=0$.
\end{definition}
\begin{definition}
$P(x,\frac{d}{dx})$ is called `Fuchsian of Type II'  
if its all singular points are 
regular singular points of 
the ODE $P(x,\frac{d}{dx})f(x)=\lambda f(x)$ for any complex number $\lambda$.
\end{definition}
These definitions do not depend on $\lambda$. From the definitions, a Fuchsin operator of type II is always  Fuchsian of Type I. The both definitions are equivalent when the coefficient function $p_M(x)$ of the highest order has no zero point whose multiplicity is greater than $M$. (However, otherwise, they are not always equivalent. For example, the differential operator $x^2\frac{d}{dx}+x$ is a Fuchsian operator of Type I but it is not of type II.) \, The definition of Fuchsian operator used in another paper~\cite{paper1} of us is 'of type I' in this paper. In the cases of Fuchsian operators of type II, 
% A differential operator 
% $P(x,\frac{d}{dx})$ is called Fuchsian
% if its all singular points are 
% regular singular points of 
% the ODE $P(x,\frac{d}{dx})f(x)=\lambda f(x)$ for $\lambda$.
% This definition does not depend on $\lambda$.
% In the Fuchsian case,
when {\bf C1$^+$} and {\bf C2$^+$} holds,
Theorem 2.7 in \cite{paper1} guarantees 
the implication  
$(iv)\Rightarrow(ii)$ in Theorem \ref{thm:equivalence1},
we obtain the following theorem:

\begin{theorem}
\Label{thm:equivalence2}
When the quintuplet
$(P(x,{\textstyle \frac{d}{dx}}),{\cal H},
\{e_n \, \}_{n=0}^{\infty},
{{\cal H}^\Diamond },
\{e_n^\Diamond \}_{n=0}^{\infty})$
satisfies 
{\bf C1-C3}, {\bf C1$^+$}, {\bf C2$^+$},and {\bf C2.1-C2.3}
and 
the differential operator 
$P(x,\frac{d}{dx})$ is Fuchsian of type II, 
then the following conditions are equivalent
for $f \in {\cal H}$:

$(ii)$: \, $f\in {\rm dom} \, B_{P,\lambda , {\cal H}, {{\cal H}^\Diamond}}$ and  $B_{P,\lambda , {\cal H}, {{\cal H}^\Diamond}}\, f = 0$

$(iii)$: \, $
f\in{\cal H} \mbox{ and } 
\{\langle f, e_n \rangle_{\cal H}\}_{n=0}^\infty \in V \cap \ell ^2(\mathbb{Z}^+)$.

$(iv)$: $f\in 
C^M(\mathbb{R}\setminus S ) 
\cap {\cal H}$ and 
$^\forall x \in \mathbb{R}\setminus S$, 
$ {P}(x,\frac{d}{dx})f(x)=\lambda f(x) $,
$S$ is the set of singular points of $(P(x,{\textstyle \frac{d}{dx}})$.
\end{theorem}

When a differential operator 
$P(x,\frac{d}{dx})$ has no singular point,
this operator is a special case of 
Fuchsian differential operators.
In this special case, 
since the set $S$ is empty, 
the relation 
$(iv)\Rightarrow(i)$ is trivial.
Then, we obtain the following theorem.

\begin{theorem}
\Label{thm:equivalence3}
When 
the quintuplet
$(P(x,{\textstyle \frac{d}{dx}}),{\cal H},
\{e_n \, \}_{n=0}^{\infty},
{{\cal H}^\Diamond },
\{e_n^\Diamond \}_{n=0}^{\infty})$
satisfies 
{\bf C1-C3}, {\bf C2$^+$},and {\bf C2.1-C2.3}
and 
the differential operator 
$P(x,\frac{d}{dx})$ is has no zero points,
then the following conditions are equivalent
for $f \in {\cal H}$:

$(i)$: \, $f\in {\rm dom} \, A_{P,{ \cal H}}$ and  
$A_{P,{ \cal H}}\, f = \lambda f$.

$(ii)$: \, $f\in {\rm dom} \, B_{P,\lambda , {\cal H}, {{\cal H}^\Diamond}}$ and  $B_{P,\lambda , {\cal H}, {{\cal H}^\Diamond}}\, f = 0$

$(iii)$: \, $
f\in{\cal H} \mbox{ and } 
\{\langle f, e_n \rangle_{\cal H}\}_{n=0}^\infty \in V \cap \ell ^2(\mathbb{Z}^+)$.

$(iv)$: $f\in 
C^M(\mathbb{R}) 
\cap {\cal H}$ and 
$^\forall x \in \mathbb{R}$, 
$ {P}(x,\frac{d}{dx})f(x)=\lambda f(x) $.
\end{theorem}

Here, we explain the structure of proof of
Theorem \ref{thm:equivalence1}.
The statements $(i)\Rightarrow(ii)$
under {\bf C1-C3} can be shown from the following lemma.

\begin{lemma}
\Label{thm:itoii}
If $f\in {\rm dom} \, A_{P,{ \cal H}}$ and  
$A_{P,{ \cal H}}\, f = \lambda f$, 
then $f\in {\rm dom} \, B_{P,\lambda , {\cal H}, {{\cal H}^\Diamond}}$ and  $B_{P,\lambda , {\cal H}, {{\cal H}^\Diamond}}\, f = 0$.
\end{lemma}

\begin{proof}
The inclusion relation ${\cal H}\subset {{\cal H}^\Diamond}$ in the sense of sets implies also that any function sequence converging for the norm $\|\cdot \|_{\cal H}$ converges for the norm $\|\cdot \|_{{{\cal H}^\Diamond}}$. Hence, from the definitions, $\, {\rm dom}\, A_{P,{ \cal H}} 
= {\rm dom}\, (A_{P,{ \cal H}}-\lambda I) \subset {\rm dom}\, B_{P,\lambda , {\cal H}, {{\cal H}^\Diamond}}$.  
 Since the equality $A_{P,{ \cal H}} f=\lambda f$ i.e. $(A_{P,{ \cal H}}-\lambda I) f=0$ implies $B_{P,\lambda , {\cal H}, {{\cal H}^\Diamond}}f=0$, this suffices for the proof of this lemma. 
\end{proof}

The statements  $(ii)\Rightarrow(iii)$ under {\bf C1-C3} 
can be shown by application of
Theorem 2.2 of ~\cite{paper1} to 
the operator $B_{P,\lambda , {\cal H}, {{\cal H}^\Diamond}}$.
The remaining part $(iii)\Rightarrow(iv)$  
will be shown in Section \ref{sec:pr-psi-in-dom-adj}
under {\bf C1}, {\bf C2}, {\bf C2$^+$} and {\bf C2.1-C2.3}.
Therefore, our remaining tasks are summarized as follows.

\begin{description}
\item[Task 1]
(Subsection \ref{subsec:s24})
Constructing
${\cal H}$, ${{\cal H}^\Diamond }$, and their CONSs satisfying Conditions  
{\bf C1}-{\bf C3}, {\bf C1$^+$}, {\bf C2$^+$}, and {\bf C2.1-C2.3}
for a differential operator $P(x,\frac{d}{dx})$ with polynomial coefficient functions $p_m(x)$.
Check of 
{\bf C1}-{\bf C2}, {\bf C1$^+$}, {\bf C2$^+$}
has been done in \cite{paper1}.

\item[Task 2]
(Section \ref{sec:sm})
Checking Condition {\bf C3} in the construction of
${\cal H}$, ${{\cal H}^\Diamond }$, and their CONSs 
given in Subsection \ref{subsec:s24}

\item[Task 3]
(Subsection \ref{subsec:s25})
Showing Theorems \ref{thm:main} and \ref{thm:main2} 
using Theorem \ref{thm:equivalence1} and 
the construction of
${\cal H}$, ${{\cal H}^\Diamond }$, and their CONSs 
given in Subsection \ref{subsec:s24}.
Indeed, when the differential operator $R(x,\frac{d}{dx})$
mentioned in Theorems \ref{thm:main} and \ref{thm:main2}
has polynomial coefficient functions,
the argument of these theorems 
are immediate from Theorem \ref{thm:equivalence1} and the above construction.
However, it is not trivial in the non-polynomial case.

\item[Task 4]
(Section \ref{sec:nenps})
Showing the relation $(iii)\Rightarrow(iv)$ mentioned in Theorem \ref{thm:equivalence1}
under {\bf C1}, {\bf C2}, {\bf C2$^+$} and {\bf C2.1-C2.3}.
\end{description}

\subsection{Function spaces and basis systems satisfying the conditions} 
\Label{subsec:s24}
In this subsection, 
we treat the case satisfying the following:
(1) The differential operator $P(x,\frac{d}{dx})$ has polynomial coefficient functions.
(2) ${\cal H}=(L_{(k_0)}^2(\mathbb{R} )$, ${\cal H}^\Diamond=L_{({k_0^\Diamond})}^2(\mathbb{R} ))$, 
(3) $k_0$ and $k_0^\Diamond$ satisfy 
${k_0^\Diamond}\le k_0-s_0$ with 
\begin{eqnarray}
s_0:=\max_m\,  (\deg {p}_m -m)\, .
\Label{eqn:def-s0}
\end{eqnarray}
The purpose of this subsection is giving CONSs of 
$L_{(k_0)}^2(\mathbb{R} )$ and $L_{({k_0^\Diamond})}^2(\mathbb{R} )$ 
such that
satisfying Conditions {\bf C1-C3}, {\bf C1$^+$}, {\bf C2$^+$}, and {\bf C2.1-C2.3}
with the above conditions.

First, we introduce 
basis systems 
$\{e_n\, |\, n\in\mathbb{Z}^+ \}$ 
and
$\{e^\Diamond_n\, |\, n\in\mathbb{Z}^+ \}$ 
of
$L_{(k_0)}^2(\mathbb{R} ) $ and $L_{({k_0^\Diamond})}^2(\mathbb{R} )$:
\begin{eqnarray}
e_n(x):= \sqrt{\frac{1}{\pi}} \,\psi_{k_0,\, \ddot{n}_{k_0,n}}(x)  ,
\quad
e^\Diamond_n(x):= \sqrt{\frac{1}{\pi}} \,\psi_{{k_0^\Diamond},\, \ddot{n}_{{k_0^\Diamond},n}}(x) 
\Label{eqn:choice-bs}
\end{eqnarray}
with
\begin{eqnarray}
\ddot{n}_{k,n}
&:=& \lfloor {\textstyle -\frac{k+1}{2}} \rfloor + (-1)^{n+k+1} \lfloor {\textstyle \frac{n+1}{2}} \rfloor
\Label{eqn:def-tilde-n} \\
\psi_{k,\, \ddot{n}}(x) 
&:=& \frac{1}{(x+i)^{k+1}} \left( \frac{x-i}{x+i} \right)^{\ddot{n}} \,\,\,\,\,\, (\ddot{n}\in\mathbb{Z} ) , 
\Label{eqn:def-psi}
\end{eqnarray}
where $\lfloor a \rfloor$ denotes the largest integer not greater than $a$.
It is easy to show that this function satisfies the following properties.
\begin{eqnarray}
\,\,\,\,\,\,\,\,\,\, \psi_{k,\, \ddot{n}}\in L_{(k)}^2(\mathbb{R} ) \, , \,\,\,   \overline{\psi_{k,\, \ddot{n}}(x)} = \psi_{k,\, -\ddot{n}-k-1}(x) \,\,\,
 {\rm and } \,\,\, \langle\psi_{k,\, \ddot{m}} \, , \, \psi_{k,\, \ddot{n}} \rangle_{(k)} = \pi \,\delta_{\ddot{m} \ddot{n}}   \,\,  . 
\Label{eqn:properties-psi}\end{eqnarray}
Moreover, they satisfy the following lemma: 
\begin{lemma}
\Label{lemma:psi-cons}
$\left\{\sqrt{\frac{1}{\pi}}\,\psi_{k,\, \ddot{n}}\, \bigl| \, \ddot{n}\in\mathbb{Z} \right\}$ is an orthonormal basis of $L_{(k)}^2(\mathbb{R} )$.
\end{lemma}
The orthonormal property is shown by (\ref{eqn:properties-psi}), though the proof of completeness is somewhat complicated. Its proof is given in Appendix A of \cite{paper1}.
This lemma guarantees {\bf C1}. 

The indices of functions in $\left\{\,\psi_{k_0,\, \ddot{n}}\, \bigl| \, \ddot{n}\in\mathbb{Z} \right\}$ are bilaterally expressed, while the indices of basis functions in  $\{e_n\, |\, n\in\mathbb{Z}^+ \}$ are unilaterally expressed, and they are 'matched' to one another by the one-to-one mapping defined by  (\ref{eqn:def-tilde-n}). In order to avoid confusion between them, in this paper, the integer indices with double dots\, $\ddot{}$\, denote the bilateral ones in $\mathbb{Z}$, in contrast to the unilateral ones (without double dots) in $\mathbb{Z}^+$. 

Since the mapping $n \to \ddot{n}_{k,n}$ is one-to-one from $\mathbb{Z}^+$ to $\mathbb{Z}$, the basis systems $\{e_n\, |\, n\in\mathbb{Z}^+ \}$ and $\{e^\Diamond_n\, |\, n\in\mathbb{Z}^+ \}$ are identical to $\left\{\sqrt{\frac{1}{\pi}}\,\psi_{k_0,\, \ddot{n}}\, \bigl| \, \ddot{n}\in\mathbb{Z} \right\}$ and  
$\left\{\sqrt{\frac{1}{\pi}}\,\psi_{{k_0^\Diamond},\, \ddot{n}}\, \bigl| \, \ddot{n}\in\mathbb{Z} \right\}$, respectively. 
%We have the following theorem:
% \begin{lemma}
% \Label{lemma:psi-cons}
% $\left\{\sqrt{\frac{1}{\pi}}\,\psi_{k,\, \ddot{n}}\, \bigl| \, \ddot{n}\in\mathbb{Z} \right\}$ is an orthonormal basis of $L_{(k)}^2(\mathbb{R} )$.
% \end{lemma}
% The orthonormal property is shown by (\ref{eqn:properties-psi}), though the proof of completeness is somewhat complicated. Its proof is given in Appendix \ref{app:proof-bs}.
% This theorem guarantees {\bf C1}.
Hence, 
%directly, 
from Lemma \ref{lemma:psi-cons},  
we have 
\begin{theorem}
\Label{thm:e-cons}
$\{e_n\, |\, n\in\mathbb{Z}^+ \}$ and $\{e^\Diamond_n\, |\, n\in\mathbb{Z}^+ \}$ are orthonormal basis systems for ${\cal H}$ and ${{\cal H}^\Diamond}$, respectively. 
\end{theorem}

 The `matched' number $\ddot{n}_{k,n}$ in (\ref{eqn:def-tilde-n}) has the property 
\begin{eqnarray}
\bigl|\, \ddot{n}_{k,n}+{\textstyle \frac{k+1}{2}}\,\, \bigr| = \left\{\begin{array}{@{\,}ll} 
\textstyle \lfloor \frac{n}{2} \rfloor +\frac{1}{2} \,\,\,\, & \,\,\,\, (k:\, {\rm even}) \\ \\ 
\textstyle \lfloor \frac{n+1}{2} \rfloor \,\,\,\, & \,\,\,\, (k:\, {\rm odd})
\end{array}\right. 
\Label{eqn:sorting-ev}\end{eqnarray}
which is used later.

As well as satisfying the orthogonality property above, they satisfy other 
orthogonality-like relations (w.r.t. other inner products) given in~\cite{paper3}, one of which is related to $\mathfrak{su}(1,1)$-number-states~\cite{SaHa}. When $k\ge 0$, as is explained in the paper~\cite{paper3} in detail, $\psi _{k,\ddot{n}}(x)$ is an `almost-sinusoidally' oscillating wavepacket with a spindle-shaped envelope $|\psi_{k,\ddot{n}}(x)| = (x^2 +1)^{-\frac{k+1}{2}} \,$,\, and its approximation to a sinusoidal wavepacket with a Gaussian envelope holds for sufficiently large $k$ with respect to the $L^2$-norm.

In the following part of this subsection, we show that 
the quintuplet
$(P,
L_{(k_0)}^2 (\mathbb{R}),
\{\sqrt{\textstyle\frac{1}{\pi}} \,\psi_{k_0,\, \ddot{n}_{k_0,n}} \}_{n=0}^\infty,
L_{(k_0^{\Diamond})}^2 (\mathbb{R}),
\{\sqrt{\textstyle\frac{1}{\pi}} \,\psi_{{k_0^\Diamond },\, \ddot{n}_{{k_0^\Diamond },n}} \}_{n=0}^\infty
)$
satisfies Conditions {\bf C2}, {\bf C2$^+$}, and {\bf C2.1-C2.3}. 
However, our proof for {\bf C3} requires several pages, 
and it will be given in Section \ref{sec:pr-psi-in-dom-adj} after the introduction of a tool for it in Section \ref{sec:sm}. 

Firstly, {\bf C2$^+$} is obvious from the definition of $\langle \cdot , \cdot\rangle _{({k_0^\Diamond})}$. Moreover, the definition of $\psi _{k,\ddot{n}}(x)$ results in the following lemma:
\begin{lemma}
\Label{lemma:psi-smooth-in-L2-k}
\hspace{3mm} $\psi _{k,\ddot{n}}\in C^\infty (\mathbb{R} ) \cap L_{(k)}^2$. 
\end{lemma} 
Since $|\psi_{{k_0^\Diamond},\ddot{n}}(x)|=(x^2+1)^{\frac{-{k_0^\Diamond}+1}{2}}$ holds for any real number $x$, {\bf C2.3} is obvious for $\tilde{a}(x)=\sqrt{\frac{1}{\pi }}\,\, (x^2+1)^{\frac{-{k_0^\Diamond}+1}{2}}$. 
In order to show {\bf C2.2}, we focus on the equality:
\begin{eqnarray}
-\frac{i}{2} \left( (x^2 +1)\, \frac{d}{dx}  + (k+1) x \, \right) \, \psi _{k,\ddot{n}} (x) = \left( \ddot{n}+\frac{k+1}{2}\right) \, \psi _{k,\ddot{n}} (x) .
\Label{eqn:char-eq}\end{eqnarray}
Then, 
the operator
$N(x,\frac{d}{dx}):=-\frac{i}{2} (x^2 +1)\, \frac{d}{dx}  + (k_0^\Diamond+1) x$
satisfies the eigen equation
$ 
N(x,\frac{d}{dx})
 \, e^\Diamond_n(x) = \lambda_n e^\Diamond_n(x)
$,
where
$\lambda _n:= \ddot{n}_{k_0^\Diamond,n} +\frac{k_0^\Diamond+1}{2}$.
Since (\ref{eqn:sorting-ev}) implies the inequality $\displaystyle |\lambda _n|> \frac{n}{2}$,
Condition {\bf C2.2} holds.
%Thus, we have shown that {\bf C2.2-C2.3} are satisfied.

Next, in order to check Conditions {\bf C2} and {\bf C2.1},
we establish some properties of $\psi_{k,\, \ddot{n}}$.
\begin{theorem}
\Label{thm:ixd-psi}
For any integer $\ddot{n}$, 
\begin{eqnarray}
\psi_{k,\, \ddot{n}}(x) 
&=& - \frac{i}{2} \left( \psi_{k-1,\, \ddot{n}} (x) - \psi_{k-1,\, \ddot{n}+1} (x), \right) 
\Label{eqn:id} \\
x\, \psi_{k,\, \ddot{n}} (x) 
&=& \frac{1}{2} \left( \psi_{k-1,\, \ddot{n}} (x) + \psi_{k-1,\, \ddot{n}+1} (x), \right)
\Label{eqn:mult} \\
{\textstyle \frac{d}{dx}}\, \psi_{k ,\, \ddot{n}} (x)   
&=& \ddot{n} \,  \psi_{k+1 ,\, \ddot{n}-1} (x) - (\ddot{n}+k+1) \, \psi_{k+1 ,\, \ddot{n}} (x) .
\Label{eqn:diff}
\end{eqnarray}
\end{theorem}
This theorem is derived directly from (\ref{eqn:def-psi}). 
A recursive use of these relations results in the following lemma:
\begin{lemma}
\Label{lemma:expansion-recursion}
Let $k_0, j, m \in \mathbb{Z}$ and ${k_0^\Diamond}\in \mathbb{Z}$. 
When $\, {k_0^\Diamond}\le k_0+m-j $, 
the function $x^j ({\textstyle \frac{d}{dx}})^m \psi_{k_0,\, \ddot{n}}(x)$ can be expressed as a linear 
combination of $\psi_{{k_0^\Diamond},\, \ddot{r}}(x)$ 
$   (\ddot{r}=\ddot{n}-m ,\, \ddot{n}-m+1,\, ...\, ,\,  \ddot{n}+m+k_0-{k_0^\Diamond})\, $ 
whose coefficients are polynomials of $n$ and $k$ with degree not greater than $m$.
\end{lemma} 

Remember that the differential operator $P(x,{\textstyle \frac{d}{dx}}) $ is given as a 
linear combination of the operators $x^j ({\textstyle \frac{d}{dx}})^m$.
By applying Lemma \ref{lemma:expansion-recursion},
(\ref{eqn:choice-bs}) -(\ref{eqn:def-psi}), we obtain in the following result:

\begin{lemma}{\rm (Theorem 4.2 (a) (b) of ~\cite{paper1})}
\Label{lemma:matrixelements}
Let $P(x,{\textstyle \frac{d}{dx}}) = \displaystyle \sum_{m=0}^M p_m(x) ({\textstyle \frac{d}{dx}})^m $. %, and let $\displaystyle\, s_0:=\max_{m} \, (\deg p_m -m) $ \, $(s_0\ge 0)$. 
When 
${p}_m(x)$ $(m=0,1,\ldots M)$ are polynomials and  
${k_0^\Diamond} \le k_0-s_0$ with $s_0$ defined as in {\rm (\ref{eqn:def-s0})},\,  $(k_0\in\mathbb{Z}^+, {k_0^\Diamond}\in\mathbb{Z} )$, the function $P(x,{\textstyle \frac{d}{dx}}) e_n(x)$ belongs to $\tilde{{\cal H}}$. 
Then, %the real number 
the complex number 
$b_m^n:=\langle Be_n,\, e^\Diamond_m\rangle =\langle P(x,{\textstyle \frac{d}{dx}}) e_n,\, e^\Diamond_m\rangle$\, $(m,n\in\mathbb{Z}^+)$ satisfies the following conditions $($a$)$ and $($b$)$:

$($a$)$ :\,  $ b_m^n=0 \mbox{ if } |m-n|>2M+k_0-{k_0^\Diamond} \,\, .$

$($b$)$ :\, There exists a polynomial $A(x)$ of degree not greater than $M$ such that \par\noindent  \hspace{1.5cm} 
$|b_m^n|\le A(n) $ for any $m,n\in\mathbb{Z}^+$.
\end{lemma} 
%The detailed proof of this is given in ~\cite{paper1}. 
Lemmata \ref{lemma:psi-cons} and \ref{lemma:matrixelements} show {\bf C2} and {\bf C2.1}. 
Thus, we have shown that 
the pair of Hilbert spaces
$(L_{(k_0)}^2(\mathbb{R} ),L_{({k_0^\Diamond})}^2(\mathbb{R} ))$
satisfies Conditions {\bf C1}, {\bf C2}, {\bf C2$^+$} and {\bf C2.1-C2.3}. 
Thus, this band-diagonal matrix $b_{m}^n$ is illustrated by Figure \ref{f1}.

\begin{figure}[htbp]
\begin{center}
\scalebox{1.0}{\includegraphics[scale=0.4]{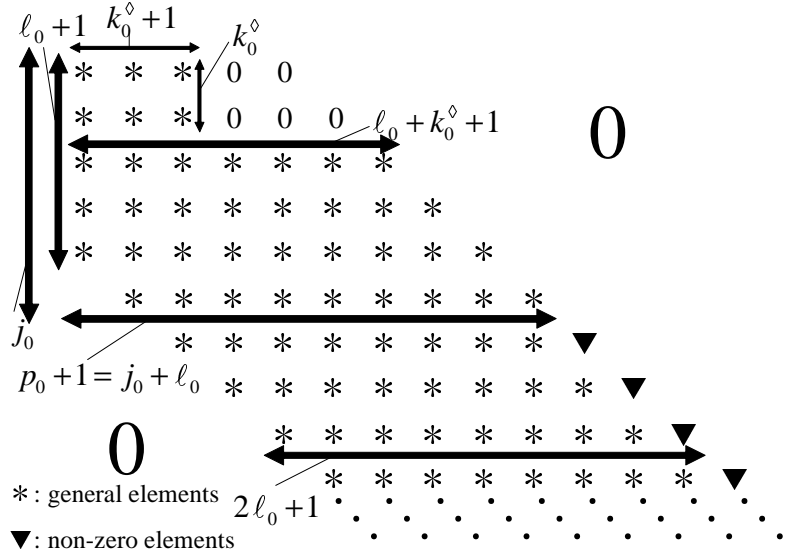}}
\end{center}
\caption{Figure of band-diagonal matrix $b_{m}^n$}
\Label{f1}
\end{figure}%

Next, we point out another property of $\psi _{k,\ddot{n}}$ related to Fourier series.
By the change of variable 
$\,x\to \theta := 2\arctan x\,$ (where $x=\tan\frac{\theta}{2}\,$), 
there is an isometric map from   
the orthonormal basis system $\{\sqrt{\frac{1}{\pi}}\,\psi_{k,\, \ddot{n}}\, | \, \ddot{n}\in\mathbb{Z} \}$ of $L_{(k)}^2(\mathbb{R} )$ to  
the orthonormal basis system of the sinusoidal waves  
$\{\frac{(-1)^n}{\sqrt{2\pi }} \,e^{in\theta} |n \in \mathbb{Z}^+ \}$ of $L^2((-\pi, \, \pi))$.
The detail of this relation is given in Appendix \ref{app:fou}. 
The same change of variable has been used for a description of analytic unit quadrature signals with nonlinear phase~\cite{Qia}~\cite{Che}, for example.
When 
%the variation of a continuous function is bounded,
a function passes Dini's test\cite{enc},  
its Fourier series satisfies point-wise convergence.
So, the above isometric correspondence between two basis systems
$\{\sqrt{\frac{1}{\pi}}\,\psi_{k,\, \ddot{n}}\, | \, \ddot{n}\in\mathbb{Z} \}$ 
and 
$\{\frac{(-1)^n}{\sqrt{2\pi }} \,e^{in\theta} |n \in \mathbb{Z}^+ \}$ 
results in Theorem \ref{thm:pointwise} in Appendix \ref{app:fou} 
which shows the point-wise convergence of the expansion of 
any once differentiable function in ${\cal H}$  
by the basis system $\{e_n\, |\, n\in\mathbb{Z}^+\}$.

\subsection{Proof of Theorems \ref{thm:main} and \ref{thm:main2}} 
\Label{subsec:s25}

Since the `if' part of Theorem \ref{thm:main} is trivial, 
it is sufficient to show the `only if' part for Theorem \ref{thm:main}.
Further, the the `only if' part for Theorem \ref{thm:main}
is a special case of Theorem \ref{thm:main2}.
Hence, we will prove only Theorem \ref{thm:main2},
which can be shown from the relation between $(i)$ and $(v)$ in the following theorem.

\begin{theorem}
Assume that a differential operator ${R}(x,{\textstyle \frac{d}{dx}})
:= \sum _{m=0}^M r_m(x) (\textstyle \frac{d}{dx})^m $ has
rational coefficient functions ${r}_m(x)$ $(m=0,1,\ldots M)$.
We denote 
the least common multiple of the denominators of ${r}_m(x)$ by $l(x)$.
For any $\lambda$, we define a differential operator ${P}(x,{\textstyle \frac{d}{dx}})
:= \sum _{m=0}^M p_m(x) (\textstyle \frac{d}{dx})^m -\lambda l(x)$ with
$p_m(x):= l(x) {r}_M(x)$.
Then, for any $k_0$, there exists an integer ${k_0^\Diamond }$ such that
the relations 
$(i)\Rightarrow(ii)\Rightarrow(iii) \Rightarrow(iv)\Rightarrow(v)$ for $f\in L_{(k_0)}^2 (\mathbb{R})$
hold.

$(i):$
$A_{R,L_{(k_0)}^2 (\mathbb{R})} f=\lambda f $.

$(ii):$
$B_{P,0, L_{(k_0)}^2 (\mathbb{R}),L_{(k_0^{\Diamond})}^2 (\mathbb{R})} f=0$.

$(iii):$
The $\ell ^2$-sequence 
$\{f_n:=\langle f, e_n\rangle_{{\cal H}}\}_{n=0}^\infty$
belongs to $V$ defined with
with the quintuplet
$(P,
L_{(k_0)}^2 (\mathbb{R}),
\{\sqrt{\textstyle\frac{1}{\pi}} \,\psi_{k_0,\, \ddot{n}_{k_0,n}} \}_{n=0}^\infty,
L_{(k_0^{\Diamond})}^2 (\mathbb{R}),
\{\sqrt{\textstyle\frac{1}{\pi}} \,\psi_{{k_0^\Diamond },\, \ddot{n}_{{k_0^\Diamond },n}} \}_{n=0}^\infty
)$.

$(iv):$
 $f \in C^M(\mathbb{R} \setminus (l\cdot{r}_M)^{-1}(0)) \cap 
L_{(k_0)}^2 (\mathbb{R})$ and 
$x \in \mathbb{R} \setminus (l\cdot{r}_M)^{-1}(0)$ satisfies
$P(x,{\textstyle \frac{d}{dx}})f(x)=0$.

$(v):$
$f \in C^M(\mathbb{R} \setminus (l\cdot{r}_M)^{-1}(0)) \cap 
L_{(k_0)}^2 (\mathbb{R})$ and 
$x \in \mathbb{R} \setminus (l\cdot{r}_M)^{-1}(0)$ satisfies
$R(x,{\textstyle \frac{d}{dx}})f(x)=\lambda f(x)$.
\end{theorem}

\begin{proof}
There exist an integer $k_1$ and a constant $c$ such that
$(x^2+1)^{k_1} (l(x))^2 \le c$.
Then, we choose $k_0^\Diamond$
with satisfying the condition 
${k_0^\Diamond }\le \min\{k_0- \max_{m} \, (\deg p_m -m), k_0+ k_1\} $.

The property 
${k_0^\Diamond }\le k_0+ k_1$ yields that
when a sequence $f_n \in L_{(k_0)}^2 (\mathbb{R})$
satisfies 
$
\|f_n\|_{(k_0)}+
\|({R}(x,{\textstyle \frac{d}{dx}})-\lambda I)  f_n\|_{(k_0^\Diamond)} \to 0$,
$
\|f_n\|_{(k_0)}+
\|{P}(x,{\textstyle \frac{d}{dx}}) f_n\|_{(k_0^\Diamond)} \to 0$.
Thus, 
the relation $(i)\Rightarrow(ii)$ holds.

Since the set of singular points of the differential operator 
${P}(x,{\textstyle \frac{d}{dx}})$ is 
$(l\cdot{r}_M)^{-1}(0)$,
the property 
${k_0^\Diamond }\le k_0- \max_{m} \, (\deg p_m -m)$
allows us to apply 
Theorem \ref{thm:equivalence1} to 
the quintuplet
$(P,
L_{(k_0)}^2 (\mathbb{R}),
\{\sqrt{\textstyle\frac{1}{\pi}} \,\psi_{k_0,\, \ddot{n}_{k_0,n}} \}_{n=0}^\infty,
L_{(k_0^{\Diamond})}^2 (\mathbb{R}),
\{\sqrt{\textstyle\frac{1}{\pi}} \,\psi_{{k_0^\Diamond },\, \ddot{n}_{{k_0^\Diamond },n}} \}_{n=0}^\infty
)$.
Then,
the relations $(ii)\Rightarrow(iii) \Rightarrow(iv)$ hold.
Since the relation $(iv)\Rightarrow(v)$ is trivial, we proved the desired arguments.

\end{proof}

\subsection{Relationship to the algorithm}
The basic framework for the proof of regularity given in this paper is the same as the framework for the algorithm proposed in~\cite{paper1} and~\cite{paper3} which yields all the solutions in 
$L_{(k_0)}^2(\mathbb{R} )$ 
% C^M(\mathbb{R} ) 
% \cap L_{(k_0)}^2(\mathbb{R} )$ 
of 
higher-order ODEs 
%the corresponding differential equation 
using only the four arithmetical operations on integers
when the ODEs have no singular points or they are Fuchsian.  
(Note that, even if the value of a solution is ambiguous at the singular points of the ODE, 
the ambiguity causes no problem for discussions in a Hilbert function space where functions are defined by equivalent classes of the quotient space by the subspace of null functions.)
% when ${p}_M(x)$ has no zero points.
% (Moreover, it yields all the solutions in $L_{(k_0)}^2(\mathbb{R} )$ of the ODE for the cases % of the Fuchsian class 
% even if it has zero points, under a modification.) \, 
This algorithm is based on the matrix representation of the operator $B_{P,\lambda , {\cal H}, {{\cal H}^\Diamond}}$
with respect to the basis systems $\{\sqrt{\frac{1}{\pi}}\,\psi_{k_0,\, \ddot{n}_{k_0,n}}\, | \, n\in\mathbb{Z}^+ \}$  and $\{\sqrt{\frac{1}{\pi}}\,\psi_{k_0^\Diamond ,\, \ddot{n}_{k_0^\Diamond,n}}\, | \, n\in\mathbb{Z}^+ \}$  under the choice of spaces  ${\cal H}=L_{(k_0)}^2(\mathbb{R} )$ and  ${\cal H}^\Diamond =L_{(k_0^\Diamond )}^2(\mathbb{R} )$ with $k_0^\Diamond \le k_0-s_0$. 
In this context, the proofs given in this paper can be interpreted as proofs of the validity of this algorithm, which guarantee the one-to-one correspondence between the square-summable vector solution of the corresponding the band-diagonal-type matrix-vector equation (simultaneous linear equations) and the true solutions in 
${\cal H}$
%$L_{(k_0)}^2(\mathbb{R} )$ 
%$C^M(\mathbb{R} )
%\cap L_{(k_0)}^2(\mathbb{R} )$ 
of the corresponding differential equation, i.e., the one-to-one correspondence between the vectors in $V \cap \ell ^2(\mathbb{Z}^+ )$ with 
$V$ defined in (\ref{eqn:sp-number-seq}) and the functions in 
$\{f\in 
C^M(\mathbb{R}\setminus S ) 
\cap {\cal H} \, | \,  ^\forall x \in \mathbb{R}\setminus S,\,  {P}(\lambda ; x,\frac{d}{dx})f(x)=0 \}$.
% $\{f \in 
% C^M(\mathbb{R} ) \cap L_{(k_0)}^2(\mathbb{R} )
% \, |\, {P}(x,{\textstyle \frac{d}{dx}})f=\lambda f\} $

Since the matrix-vector equation $\displaystyle\sum_n b_m^n f_n =0$ corresponding to the ODE is infinite-dimensional, we should be careful of 
whether or not the vector corresponding to any solution in  $L_{(k_0)}^2(\mathbb{R} )$ of the ODE  ${P}(x,{\textstyle \frac{d}{dx}})f=\lambda f$ (on $\mathbb{R}\setminus S$) 
always satisfies the matrix-vector equation $\displaystyle\sum_n b_m^n f_n =0$. 
In this context, the proofs of $(iv)\Longrightarrow (ii)\Longrightarrow (iii)$ 
can be regarded as the proof of the validity of the matrix-vector representation of the ODE. 

On the other hand, 
in infinite-dimensional case, all the solutions of the matrix-vector equation do not necessarily  correspond to the true solutions of the ODE. Actually, as is shown in~\cite{paper1},  
there are vectors in $V$ which do not correspond to any true solution in 
$L_{(k_0)}^2(\mathbb{R} )$ of the ODE  
(on $\mathbb{R}\setminus S$); nevertheless there is no such vector in 
$V\cap \ell ^2(\mathbb{Z}^+ )$. However, in~\cite{paper1}, the statement $(iii)\Longrightarrow (iv)$ is  assumed only  as a condition, which is {\bf C4} of ~\cite{paper1}, and its proof is omitted in that paper.
In this context, 
the proof of $(iii)\Longrightarrow (iv)$ can be regarded as a proof of the non-existence of extra solutions in $L_{(k_0)}^2(\mathbb{R} ) $ in our method which do not correspond to any 
solution in  
$L_{(k_0)}^2(\mathbb{R} )$ of the corresponding ODE  
(on $\mathbb{R}\setminus S$). 
% There are vectors in $V$ which do not correspond to any true solution in 
% $L_{(k_0)}^2(\mathbb{R} )$ of the differential equation 
% \R{
% (on $\mathbb{R}\setminus S$) 
% }
% as is shown in~\cite{paper1}; nevertheless there is no such vector in 
% $V\cap \ell ^2(\mathbb{Z}^+ 
In the proposed algorithm, we utilizes a method for the removal of the non-square-summable components from the vectors in $V$, and hence we can obtain approximations for only the true solutions $L_{(k_0)}^2(\mathbb{R} )$ of the differential equation with high accuracy. 
%In~\cite{paper1}, the statement $(iii)\Longrightarrow (iv)$ is  assumed only  as a condition, which is {\bf C4} of ~\cite{paper1} and whose proof is omitted in that paper.

Thus, the proofs in this paper guarantee also the one-to-one correspondence between the functions obtained by this integer-type algorithm and the true solutions in $L_{(k_0)}^2(\mathbb{R} )$ of the differential equation
(on $\mathbb{R}\setminus S$) 
. From this point of view, this paper 
contains the proofs of some propositions required in~\cite{paper1}, which was omitted there. They will be given in Subsections \ref{sec:pr-psi-in-dom-adj} of this paper. 
%is useful to show the validity of~\cite{paper1}, because three proofs omitted there will be given in 
%Subsections \ref{sec:pr-psi-in-dom-adj} of this paper. 

% Even when the coefficient polynomial of the highest order has zero points, as is shown in Appendix \ref{sec:fuchs}, 
% if the ODE belongs to the Fuchsian class, the following modification of the usage of the algorithm guarantees the one-to-one correspondence for $L^2$-solutions similarly to the above. That is, the one-to-one correspondence is guaranteed when we apply the algorithm to a Fuchsian-type ODE after multiplying a polynomial from the left. 
% Note that, even if the value of a solution is ambiguous at the zero points of the coefficient polynomial, it causes no problem for discussions in a Hilbert function space where functions are defined by equivalent classes of the quotient space by the subspace of null functions.
% For non-Fuchsian cases where 
% ${p}_M(x)$ has zero points, 
More generally, for non-Fuchsian cases, we can use at least 
%instead of the one-to-one correspondence, 
% between the vectors in $V\cap \ell ^2(\mathbb{Z}^+ )$ and the functions in $\{f \in C^M(\mathbb{R} ) \cap L_{(k_0)}^2(\mathbb{R} )\, |\, {P}(x,{\textstyle \frac{d}{dx}})f=\lambda f\} $, 
%we use 
the statement  $(iii)\Rightarrow(iv )$ of Theorem \ref{thm:equivalence1}. This statement guarantees that all the solution obtained by this algorithm approximately coincide with true solutions of the differential equation in any open interval between two adjacent 
singular points of the ODE, 
%zero points of ${p}_M(x)$, 
and it guarantees that this algorithm yields at least all the solutions in  
$C(\mathbb{R}\setminus S) \cap ({\rm dom } \, B_{P,\lambda, { \cal H}, {\cal H}^\Diamond})$ 
%${C}^M(\mathbb{R}\backslash {p}_M^{-1}(0)) \cap ({\rm dom } \, A_{P,{ \cal H}})$ 
of the ODE 
(on $\mathbb{R}\setminus S$).

\section{A `kind of smoothing operator' and Condition C3}
\Label{sec:sm}
\subsection{A `kind of smoothing operator' for blurring endpoints}

In order to show {\bf C3}, %whether a function belongs to the domain of the adjoint operator of a operator defined as the action of $P(x,\frac{d}{dx})$, 
we have to check whether the contribution of the difference terms between two endpoints in the `integration by parts' vanish or not as the endpoints tend to $\pm\infty$. Usually, for functions in a Hilbert space in general, it is difficult to show this vanishing by a direct method because the normalizability does not always imply smooth decays for large $|x|$ but may possibly allow long-lasting sparse oscillations with undesired peak amplitudes. For the proof based upon this vanishing, here we will introduce a convenient operator $T$ which `blurs' the two endpoints.
\begin{definition}
\Label{definition:S}
On a space in general of locally integrable functions, define the linear operator $T$ by 
\begin{eqnarray*}
\left( T f \right) (x) := \left\{\begin{array}{@{\,}ll}\displaystyle 
\,\, \frac{1}{x} \int_{x}^{2x} f(u) \, du \,\, & ({\rm if} \,\,\, x\ne 0) \\ \\ \displaystyle 
\,\, f(0) & ({\rm if} \,\,\, x=0) \,\,\, . 
\end{array}\right. 
\end{eqnarray*}
\end{definition}
\begin{lemma}
\Label{lemma:S}
The operator $T$ defined above satisfies the following properties: 

$\,\left( T f \right) (x)$ is $(m+1)$-times continuously differentiable in $\mathbb{R} \backslash \{ 0\}$ if $f(x)$ is at least $m$-times continuously differentiable in $\mathbb{R} \backslash \{ 0\}$.   Moreover, 

\begin{eqnarray}
\left(Tg\right)(x) &=& \left(Tf\right)(cx) 
\,\,\,\,\,\,\, {\it if} \,\,\,\,\,\,\,\,
g(x)=f(cx) \,\,\,\,\,\, (c:\, {\it nonzero \,\, real \,\, constant}) ,
\Label{eqn:sc_sm} \\
\hspace{5mm}\lim_{x\to\pm\infty } \left(Tf\right)(x) 
&=& \lim_{x\to\pm\infty } f(x) 
\,\,\,\,\,\,\, {\it if} \,\,\,\,\,\,\,\,
 ^\exists \lim_{x\to\pm\infty } f(x)  , 
\Label{eqn:limit_sm} \\
\left| \, \left(Tf\right)(x) \, \right|
& \le & \left(T|f|\,\right)(x) \le  \left(T|g|\,\right)(x) 
\Label{eqn:abs_sm}
\\
&&\hspace{10ex}
{\it if} \,\,\,
 \,\,\, |f(u)| \le |g(u)| \,\,\, {\it holds \,\, for} \,\, |x| \le \,  |u| \le 2|x| \,  .  \nonumber
\end{eqnarray}
\end{lemma} 

Here we omit a discussion about differentiability at $x=0$, which has nothing to do with the proofs in this paper. The proof of this lemma is derived directly from the definition of $T$, where the negative sign cancels out when $x<0$ because then $x>2x$. 
The property (\ref{eqn:limit_sm}) in Lemma \ref{lemma:S} is very important for our purpose because it results in the following lemma:
\begin{lemma}
\Label{lemma:endpoints-vanish}
$($In the following, $f^{(n)}$ denotes $(\frac{d}{dx})^nf$ for $n\in \mathbb{Z}^+$.$)$  Let $m\in\mathbb{Z}^+$. For functions $f, \, g \in C^m(\mathbb{R})$, if there exist nonnegative integers $n_r$ $(r=0,1,2,...,m-1)$ such that  
\begin{eqnarray*}
\lim_{x\to\pm\infty}\bigl(T^{n_r}(f^{(r)} g^{(m-r-1)})\bigr)(x)=0 \,\,  
\Bigl(\, {\rm with} \, \bigl( f^{(r)} g^{(m-r-1)}\bigr)(x):=f^{(r)}(x)\, g^{(m-r-1)}(x) \Bigr) \,  
\end{eqnarray*}
for $r=0,1,2,...,m-1$ and both of $\displaystyle \int_{-\infty}^\infty f(x)\, g^{(m)}(x)  dx $ and $\displaystyle \int_{-\infty}^\infty f^{(m)}(x) \, g(x)  dx $ exist, 
then
\begin{eqnarray*}
\int_{-\infty}^\infty f(x)\,\, g^{(m)}(x) \, dx = (-1)^m \int_{-\infty}^\infty f^{(m)}(x) \,\, g(x) \, dx  \, .
\end{eqnarray*}
\end{lemma}  
{\em Proof of Lemma }\ref{lemma:endpoints-vanish}: \quad
\rm
Define 
\begin{eqnarray*}
\tilde{Y}(x):=\int_{-x}^x f(u)\,\, g^{(m)}(u) \, du 
\,\,\,\,\,\,\,\,\,\,\,\,\, {\rm and} \,\,\,\,\,\,\,\,\,\,\,\,\,
\tilde{Z}(x):=\int_{-x}^x f^{(m)}(u)\,\, g(u) \, du \,\, .
\end{eqnarray*}
Then, integrating by parts (which is always applicable to integrations over a finite interval $[-x,x]\, $), 
\begin{eqnarray*}
\tilde{W}(x)&:=& \tilde{Y}(x)-(-1)^m \tilde{Z}(x)\\
  &=& \sum_{r=0}^{m-1} (-1)^r \Bigl(f^{(r)}(x)\,  g^{(m-r-1)}(x) - f^{(r)}(-x)\,  g^{(m-r-1)}(-x)\Bigr). 
\end{eqnarray*}
Since a recursive use of (\ref{eqn:limit_sm}) in Lemma \ref{lemma:S} results in 
\begin{eqnarray*}
\lim_{x\to\infty}\bigl(T^{n}f\bigr)(x)=0 \,\,\,\,\,\,\,\,\,\,\, {\rm if} \,\,\,\,\,\,\, ^\exists\ell \in \{0,1,2,...,n-1\}  \,\,\,\, {\rm s.t.} \, \lim_{x\to\pm\infty}\bigl(T^{\ell }f\bigr)(x)=0  \,\, ,
\end{eqnarray*}
with $\displaystyle n:=\max_r n_r$, we have $\, \displaystyle \lim_{x\to\infty}\bigl(T^{n}(f^{(r)} g^{(m-r-1))})\bigr)(\pm x)=0 \,\,\, {\rm for} \,\,\, r=0,1,...,m-1 \,$. Hence $\displaystyle\lim_{x\to\infty}\bigl(T^{n}\tilde{W}\bigr)(x)=0$. On the other hand, 
\begin{eqnarray*}
\lim_{x\to\infty}\bigl(T^{n}\tilde{Y}\bigr)(x)
&=& \lim_{x\to\infty} \tilde{Y}(x)= \int_{-\infty}^\infty f(x)\,\, g^{(m)}(x) \, dx \\  
\lim_{x\to\infty}\bigl(T^{n}\tilde{Z}\bigr)(x)
&=& \lim_{x\to\infty} \tilde{Z}(x)=\int_{-\infty}^\infty f^{(m)}(x)\,\, g(x) \, dx . 
\end{eqnarray*}
From these facts, $\displaystyle\lim_{x\to\infty}\bigl(T^{n}\tilde{W}\bigr)(x)=0$ results in the conclusion of the lemma, because $T^n$ is linear.
\hfill\endproof

There are some other properties of $T$, useful for the proofs, which are summarized in the following lemmata:
\begin{lemma}
\Label{lemma:conv-sml2}
Let $k\in\mathbb{Z}$. For any locally integrable $f$ in $L_{(k)}^2(\mathbb{R} )$, with $\, p(x):=x^k f(x)\,$,  $\displaystyle\,\lim_{x\to\pm\infty } \left(Tp\right)(x) = \lim_{x\to\pm\infty } \left(T\, |p|\,\right)(x) =0\,$.
\end{lemma} 
{\em Proof of Lemma }\ref{lemma:conv-sml2}: \quad
\rm 
From the Schwartz inequality, for $\, x\ne 0\,$, 
\begin{eqnarray*}
 \bigl|\left(T\, |p|\,\right)(x)\,\bigr| 
&=&\frac{1}{|x|} \,\left|\int_x^{2x} |u^k f(u)|\, du\right| 
\\
& \le & \frac{1}{|x|} \, \sqrt{|x| \cdot \left|\int_x^{2x}|u^k f(u)|^2\, du\right|} 
=  \sqrt{\frac{1}{|x|} \,\left|\int_x^{2x} u^{2k}|f(u)|^2\, du\right|} 
\end{eqnarray*}
Let $\displaystyle\, C:=\int_{-\infty}^\infty |f(u)|^2\,(u^2+1)^k \,  du  \,$. (If $\,f\in L_{(k)}^2(\mathbb{R} )\,$, $C$ should be finite.)\, Then 
\begin{eqnarray*}
\left|\int_x^{2x} u^{2k} |f(u)|^2 \, du\right|
&\le & \left(\max \bigl(\,1,\, ({\textstyle\frac{x^2}{x^2+1}})^k \bigr)\right)\,\, \left|\int_x^{2x} |f(u)|^2 \, (u^2+1)^k \, du \right| \\ 
& \le & C \,\max \bigl(\,1,\, ({\textstyle\frac{x^2}{x^2+1}})^k \bigr).
\end{eqnarray*}
Hence, if $\,f\in L_{(k)}^2(\mathbb{R} )\,$, then 
$\displaystyle\bigl|\left(T\, |p|\,\right)(x)\bigr| \le \frac{C\,\max \bigl(\,1,\,\, ({\textstyle\frac{x^2}{x^2+1}})^k \bigr)}{\sqrt{|x|}}\,$ holds for $\,x\ne 0$. Since $\displaystyle \lim _{x\pm\infty }\frac{(\frac{x^2}{x^2+1})^k}{\sqrt{|x|}}=0$ for any $k\in \mathbb{Z}$, with (\ref{eqn:abs_sm}), the proof is complete.  
\hfill\endproof
\par 
\begin{lemma}
\Label{lemma:prelim-inv-lhopital}
For $m\in\mathbb{Z}^+$, if $\, f\in C^1(\mathbb{R} )$ satisfies $\displaystyle\,\lim_{x\to\pm\infty } \left( T^m f \right) (x) = 0\,$, then \par\noindent  $\displaystyle\, \lim_{x\to\pm\infty } \left(T^{m+1} g\right)(x) = 0 \,$ for $\, g(x):=x{\textstyle\frac{d}{dx}}\, f(x)\, $.   
\end{lemma} 
{\em Proof of Lemma }\ref{lemma:prelim-inv-lhopital}: \quad
\rm 
Since 
\begin{eqnarray*}
 \left(Tg\right)(x) 
&=& \frac{1}{x} \int_{x}^{2x} u {\textstyle \frac{d}{du}} f(u) \, du 
\\ 
&= &\frac{1}{x} \left( \, (2x)f(2x)-xf(x) - \int_{x}^{2x} \, f(u) \, du \, \right) 
= 2  f(2x) - f(x) - \bigl(T f\bigr)(x)  , 
\end{eqnarray*}
from Definition \ref{definition:S}, (\ref{eqn:sc_sm}) and (\ref{eqn:limit_sm}), we have 
\begin{eqnarray*}
 \lim_{x\to\pm\infty } \left(T^{m+1} g\right)(x) 
&=& \lim_{x\to\pm\infty } \left(T^mT g\right)(x) 
\\ 
&= &\lim_{x\to\pm\infty } \left[\, 2 \left( \, T^m f \right) (2x) - \left( \, T^m f \right) (x) - \left(T^{m+1} f\bigr)(x) \right)\,\right] = 0 .
\end{eqnarray*}
\hfill\endproof
\par

\begin{lemma} 
\Label{lemma:inv-lhopital}
Let $f\in C^M$ be a locally integrable function satisfying
$\displaystyle \lim_{x\to\pm\infty } \left(T\, |p \, |\, \right)(x) = 0 \, $ with $\, p(x):= x^k f(x)\, $, and let $g(x)$ satisfy the following conditions {\rm ($\alpha$)-($\gamma$)}: \par\noindent  
${}\,\,\,\, {\rm (\alpha )}$ There exists $x_0 >0$ such that $g(x)$ is at least once continuously differentiable \par\noindent \hspace{8mm} for all $x$ such that $|x|> x_0$ \par\noindent   ${}\,\,\,\,{\rm (\beta )}\displaystyle\, \limsup_{x\to\pm\infty } |g(x)| < \infty \, $ \par\noindent   ${}\,\,\,\, {\rm (\gamma )}\displaystyle\, \limsup_{x\to\pm\infty } |x {\textstyle\frac{d}{dx}}\, g(x)| < \infty \, $.  \par\noindent  Then, for the functions $h_n(x):= x^{n+k} g(x) \,  ({\textstyle\frac{d}{dx}})^n f(x)\,$ $(n\in\mathbb{Z}^+ )$, the convergence $\displaystyle\, \lim_{x\to\pm\infty } \left(T^{n+1} h_n \right)(x) = 0 \,$ holds. 
\end{lemma} 
{\em Proof of Lemma }\ref{lemma:inv-lhopital}: \quad 
\rm 
The proof is by mathematical induction. 

Firstly, for the case with $n=0$ $\displaystyle\bigl({\rm where}\,\,  h_0(x)=x^k g(x)f(x)\, \bigr)$, from (\ref{eqn:abs_sm}), the theorem of the lemma holds, because the conditions of the lemma guarantee that 
\begin{eqnarray*}
 ^\exists x_c >0 \,\,\, {\rm and } \,\,\, ^\exists C>0 && \\ \,\,\,\, {\rm s.t. } && \,\,\,\, ^\forall x<-x_c \,\,\, {\rm and } \,\,\, ^\forall x>x_c \, , \,\, |g(x)|<C \,\,\,\, {\rm i.e.} \,\,\, |h_0 (x)|<C|\, x^k f(x)|  . 
\end{eqnarray*}

Next, assume that the theorem of the lemma holds for $n=0,1,2, ...,n'$. The following discussion refers only to values of $x$ such that $|x|> x_0$ where $g(x)$ is differentiable, which creates no problem for statements about the limit as $x\to\pm\infty$.  From this assumption and Lemma \ref{lemma:prelim-inv-lhopital}, 
$\displaystyle\,\lim_{x\to\pm\infty } \left(T^{\, n'+2} b_{n'} \right)(x) = 0\,$ with  $\, b_{n'}(x):=x{\textstyle\frac{d}{dx}}\, h_{n'}(x)\, $. Here, let 
\begin{eqnarray*}
q(x):=x^{n'+k} \Bigl( \, (n'+k) g(x) + x {\textstyle \frac{d}{dx}} \, g(x) \Bigl) \cdot \Bigl( ({\textstyle \frac{d}{dx}})^{n'}f(x)\,\Bigr) \,\, .
\end{eqnarray*}
Then, since 
\begin{eqnarray*}
\left(b_{n'} \right)(x) 
&=& \Bigl(\, x^{n'+k+1} \, g(x) \,\,  ({\textstyle \frac{d}{dx}})^{n'+1} f(x)\,\Bigr) 
+ \Bigl(\, x\, {\textstyle \frac{d}{dx}} \bigl(x^{n'+k} g(x)\bigr) \,\Bigr) \cdot \Bigl( ({\textstyle \frac{d}{dx}})^{n'}f(x)\,\Bigr) 
\\ 
&=& h_{n'+1} (x)+ q(x) ,
\end{eqnarray*}
we obtain 
\begin{eqnarray*}
\lim_{x\to\pm\infty } \left[ \left(T^{\, n'+2} h_{n'+1} \, \right)(x) 
+ \left(T^{\, n'+2} q \, \right)(x) \, \right] = 0 . 
\end{eqnarray*}
Since the trigonometric inequality and the conditions of the lemma imply that $\displaystyle\, \limsup_{x\to\pm\infty } |(n'+k)g(x) + x {\textstyle\frac{d}{dx}}\, g(x)|<\infty $, the statement of this lemma with $n=n'$ and (\ref{eqn:abs_sm}) result in $\displaystyle\, \lim_{x\to\pm\infty } \left(T^{\, n'+1} q \right)(x) = 0 \, $, and hence $\displaystyle\, \lim_{x\to\pm\infty } \left(T^{\, n'+2} q \right)(x) = 0 \,$ by (\ref{eqn:limit_sm}), From these relations, $\displaystyle\, \lim_{x\to\pm\infty } \left(T^{\, n'+2} h_{n'+1} \right)(x) = 0 \,$ i.e. the statement of the lemma holds for $n=n'+1$. 
\hfill\endproof
\par

\subsection{Proof of Condition {\bf C3}}
\Label{sec:pr-psi-in-dom-adj}

In this section, we will prove the following theorem, which shows {\bf C3}:
\begin{theorem}
{\rm (Theorem 4.8 of ~\cite{paper1})} \, 
\Label{thm:basis-in-dom-adj}
Let $P(x,{\textstyle \frac{d}{dx}}) = \displaystyle \sum_{m=0}^M p_m(x) ({\textstyle \frac{d}{dx}})^m $ with $\displaystyle p_m(x):=\sum_{j=0}^{\deg p_m} p_{m,j} \, x^j$, and let  $s_1\ge s_0$ $\bigl($with $s_0$ defined in {\rm (\ref{eqn:def-s0}}$)\bigr)$. %$\displaystyle\, s:=\max_{m} \, (\deg p_m -m) $ \, $(s\ge 0)$. %When $k\ge s$, then, 
Then, for the closed extension $B$ with respect to the graph norm of the operator $\tilde{B}$ defined by the action of $P(x,{\textstyle \frac{d}{dx}})$ 
with domain 
\begin{eqnarray*}
D(\tilde{B}) =\{f\in C^M \cap L_{(k_0)}^2(\mathbb{R} ) \,\, |\,\, \tilde{B}f\in 
%L_{(k-s_1)}^2 (\mathbb{R}) 
L_{(k_0-s_1)}^2 (\mathbb{R})   
\} \,\, , 
\end{eqnarray*}
and for the closed extension $C$ with respect to the graph norm of the operator $\tilde{C}$ defined by 
\begin{eqnarray*}
\left (\tilde{C} g \right)(x) := \sum_{m=0}^M \sum_{j=0}^{\deg p_m} \, (-1)^m \overline{p_{m,j}} \, (x^2+1)^{-k_0} ({\textstyle\frac{d}{dx}})^m \Bigl( x^j (x^2+1)^{k_0-s_1} \, g(x) \, \Bigr)
\end{eqnarray*}
with domain 
\begin{eqnarray*}
D(\tilde{C}) =\{f\in C^M \cap 
%L_{(k-s_1)}^2(\mathbb{R} 
L_{(k_0-s_1)}^2 (\mathbb{R})   
) \,\, |\,\, \tilde{C}f\in L_{(k_0)}^2 (\mathbb{R}) \} \,\, , 
\end{eqnarray*}
the following holds: 
\begin{eqnarray*}
 ^\forall f \in {\rm dom}\, \tilde{B} \,\,\, {\it and } \,\,\, ^\forall n \in \mathbb{Z}\,  , \,\,\,\,\, 
\bigl( \, B \, f , \, \psi_{k_0-s_1,\, \ddot{n}} \, \bigr)
%_{(k-s_1)} 
_{(k_0-s_1)}  
= \Bigl( \, f, \, C \, \psi_{k_0-s_1,\, \ddot{n}} \, \Bigr)_{(k_0)} \, .
\end{eqnarray*}
\end{theorem}
This theorem (together with results on limits of function sequences) implies that the basis functions of 
${\cal H}^\Diamond $ 
belong to the domain of the adjoint of $B$
under the above choices of function spaces and basis systems. 
This theorem is essential in order to show 
the statement $(ii)\Longrightarrow(iii)$ of Theorems \ref{thm:equivalence1},
and it guarantees that the corresponding number sequence $\{f_n\}_{n=0}^\infty$ of any true solution $f$ in 
${C}^M(\mathbb{R}\backslash {p}_M^{-1}(0)) \cap ({\rm dom } \, A_{P,{\cal H}})$
(or in $C^M(\mathbb{R} )\cap {\cal H}$ for the cases where ${p}_M(x)$ has no zero point) 
%$C^M(\mathbb{R} )\cap {\cal H}$
of the differential equation always satisfies 
the simultaneous linear equations $\sum_n b_m^n f_n =0$ $(m\in\mathbb{Z}^+)$.

Before the proof, we establish the following preliminary lemma:
\begin{lemma}
\Label{lemma:rational-lambda}
Let $k,\ddot{n}\in\mathbb{Z}$ and $j, m \in \mathbb{Z}^+$, and define 
%\par\noindent  
$\nu _{k\,\ddot{n}}:=\max (\ddot{n}+k+1,\, -\ddot{n}, \, k+1)$. Then, for the function  %\par\noindent  
$\displaystyle \, \lambda_{j,k,\ddot{n}}^{(m)}(x) :=({\textstyle\frac{d}{dx}})^m \left( x^j (x^2+1)^k \overline{\psi_{k,\, \ddot{n}}(x)}\right)\, $, the function 
%\par\noindent  
$R_{m,j,k,\ddot{n}}(x):=(x^2+1)^{\nu _{k\,\ddot{n}}+m} \, \lambda_{j,k,\ddot{n}}^{(m)}(x)$ is a polynomial in $x$ and its degree is not greater than $2\nu _{k\,\ddot{n}}+m+j+k-1$.
\end{lemma} 
{\em Proof of Lemma }\ref{lemma:rational-lambda}: \quad
\rm
From the definition (\ref{eqn:def-psi}) of $\psi _{k,\ddot{n}}(x)$, the function $(x^2+1)^{\nu _{k\,\ddot{n}}} \overline{\psi _{k,\ddot{n}}(x)}$ is a polynomial in $x$ and its degree is $2\nu _{k\,\ddot{n}} -k-1$, because the degrees of the factors $(x\pm i)$ in the denominator of $\psi _{k,\ddot{n}}(x)$ are not greater than $\nu _{k\,\ddot{n}}$ and the difference between the degree of the numerator and that of the denominator of $\psi _{k,\ddot{n}}(x)$ is $k+1$. Hence, the function  $T_{m,j,k,\ddot{n}}(x):=({\textstyle\frac{d}{dx}})^m x^j (x^2+1)^{\nu _{k\,\ddot{n}}} \overline{\psi_{k,\, \ddot{n}}(x)}$ is a polynomial in $x$ and its degree is $2\nu _{k\,\ddot{n}}-m+j-k-1$ when  $m\le 2\nu _{k\,\ddot{n}}+j-k-1$, while $T_{m,j,k,\ddot{n}}(x)=0$ when $m> 2\nu _{k\,\ddot{n}}+j-k-1$.

On the other hand, the function 
$\displaystyle T_{m,k,\ddot{n}}(x):=(x^2+1)^{m+k} \overline{\psi_{k,\, \ddot{n}}(x)} \, \left(({\textstyle\frac{d}{dx}})^m (x^2+1 )^{\nu _{k\,\ddot{n}}-k}\right)$ is a polynomial in $x$, because $({\textstyle\frac{d}{dx}})^m (x^2+1 )^{\nu _{k\,\ddot{n}}-k}$ contains the factor $(x^2+1)^{\nu _{k\,\ddot{n}}-m-k}$ when $m\le \nu _{k\,\ddot{n}}-k$. Here, when $m\le \nu _{k\,\ddot{n}}-k$, the degree of $(x^2+1)^{-\nu _{k\,\ddot{n}}+m+k}({\textstyle\frac{d}{dx}})^m (x^2+1 )^{\nu _{k\,\ddot{n}}-k}$ (which is a polynomial) is $m$. When $m>\nu _{k\,\ddot{n}}-k$, the degrees of $(x^2+1)^{m+k}\overline{\psi_{k,\, \ddot{n}}(x)}$ (which is a polynomial) and $({\textstyle\frac{d}{dx}})^m (x^2+1 )^{\nu _{k\,\ddot{n}}-k}$ are $2m+k-1$ and $2\nu _{k\,\ddot{n}}-2k -m$, respectively. From these facts, we can easily show that the degree of $T_{m,k,\ddot{n}}(x)$ is $2\nu _{k\,\ddot{n}}+m-k-1$,

Since $R_{m,j,k,\ddot{n}}(x)=(x^2+1)^{k+m}T_{m,j,k,\ddot{n}}(x)-x^jT_{m.k.\ddot{n}}(x)$, the calculations of the degrees of polynomials 
\begin{eqnarray*}
2(m+k)+(2\nu _{k\,\ddot{n}}-m+j-k-1)&=&j+ (2\nu _{k\,\ddot{n}} -k-1)+m \\ =j+(2\nu _{k\,\ddot{n}} +m-k-1) &= &2\nu _{k\,\ddot{n}}+m+j+k-1
\end{eqnarray*}
lead us to the statement of the lemma.
\hfill\endproof
\par

By means of the lemmata in Section \ref{sec:sm} about the operator $T$ and the above Lemma \ref{lemma:rational-lambda}, the proof of Lemma 4.8 of the paper~\cite{paper1} is constructed as follows: 

{\em Proof of Theorem }\ref{thm:basis-in-dom-adj}: \quad
\rm 
For $\lambda_{j,\, k_0-s_1,\, \ddot{n}}^{(m)}(x) :=({\textstyle\frac{d}{dx}})^m \left( x^j (x^2+1)^{k_0-s_1} \overline{\psi_{k_0,\, \ddot{n}}(x)}\right)\, $,  Lemma \ref{lemma:rational-lambda} implies that there exist finite $K,\xi >0$ such that $\displaystyle \bigl|\, \lambda_{j,\,k_0-s_1,\, \ddot{n}}^{(m)}(x)\, \bigr| \le K (\sqrt{x^2+1}\, )^{k_0-s_1-m+j-1}$ for $|x|>\xi $ \,\,  i.e. 
$\displaystyle | \, (x^2+1)^{-k_0} ({\textstyle\frac{d}{dx}})^m \Bigl( x^j (x^2+1)^{k_0-s_1} \psi _{k_0-s_1,\, \ddot{n}} \Bigr) |
 \le \frac{K}{(\sqrt{x^2+1}\, )^{k_0+s_1+m-j+1}}$ for $|x|>\xi $,  
Hence, there exists a real number $K'$ such that 
\begin{eqnarray*}
&&\left| \sum_{m=0}^M \sum_{j=0}^{\deg p_m} (-1)^m \overline{p_{m,j}}  \, (x^2+1)^{-k_0} ({\textstyle\frac{d}{dx}})^m \Bigl( x^j (x^2+1)^{k_0-s_1} \psi _{k_0-s_1,\, \ddot{n}} \Bigr)(x)\right|  \\ 
&\le &
\sum_{m=0}^M \sum_{j=0}^{\deg p_m} |p_{m,j}| \cdot \left| \, (x^2+1)^{-k_0} ({\textstyle\frac{d}{dx}})^m \Bigl( x^j (x^2+1)^{k_0-s_1} \psi _{k_0-s_1,\, \ddot{n}} \Bigr)(x) \right| \\  
&\le &
\sum_{j=0}^{\deg p_m} |p_{m,j}| \cdot  \frac{K'}{(\sqrt{x^2+1}\, )^{k_0+s_1+m-j+1}} \,\,\,\,\,\,\,\,\,\, \mbox{ for $|x|>\xi $.}
\end{eqnarray*}
Since  $s_1+m-j\ge s_0+m-\deg p_m \ge 0$ is satisfied for $j\le \deg p_m$, it is easily shown that 
\begin{eqnarray*}
\int_{-\infty}^\infty \left| \sum_{m=0}^M \sum_{j=0}^{\deg p_m} (-1)^m \overline{p_{m,j}}\,   (x^2+1)^{-k_0} ({\textstyle\frac{d}{dx}})^m \Bigl( x^j (x^2+1)^{k_0-s_1} \psi _{k_0-s_1,\, \ddot{n}} \Bigr)(x)\right|^2 && \\ 
\cdot \, (x^2+1)^{k_0} dx &<& \infty   
\end{eqnarray*}
from the above inequality, i.e., $\psi _{k_0-s_1,\, \ddot{n}}\in D(\tilde{C})$ .
Hence, $\tilde{C}\psi _{k_0-s_1,\, \ddot{n}}$ is well defined and 
\begin{eqnarray*}
\bigl(\tilde{C}\psi _{k_0-s_1,\, \ddot{n}}\bigr)(x)  = \sum_{m=0}^M \sum_{j=0}^{\deg p_m} (-1)^m \, \overline{p_{m,j}}\,\, (x^2+1)^{-k_0} \,\overline{\lambda_{j,\,k_0-s_1,\, \ddot{n}}^{(m)}(x)} \,\, .
\end{eqnarray*}

(In the following, the suffixes for $j$, $m$, $k_0$, $s_1$ and $\ddot{n}$ are often omitted if unnecessary for simplicity.) \,  

Let $f\in D(\tilde{B})\,$. Then, for 
\begin{eqnarray*}
Z(x) := \int_{-x}^x \bigl(\tilde{B}f\bigr)(u) \,\, \overline{\psi_{k_0-j,\, \ddot{n}}(u)} \,\, (u^2+1)^{k_0-j}\, du \,\, ,  
\end{eqnarray*} 
the convergence $\, \displaystyle \lim_{x\to\infty } Z(x) = \bigl( \, \tilde{B}f, \, \psi_{k_0-j,\, \ddot{n}}\, \bigr)_{(k_0-j)}$ holds because
$\tilde{B}f\in \tilde{{\cal H}}=L_{(k_0-s_1)}(\mathbb{R} )\, $ and $\, \psi _{k_0-s_1,\, \ddot{n}}\in L_{(k_0-s_1)}(\mathbb{R} )$. Next, define 
\begin{eqnarray*}
Y(x)  &:=& \int_{-x}^x  f(u) \,\, \overline{\bigl(\tilde{C}\psi _{k_0-s_1,\, \ddot{n}} \bigr) (u)}  \,\,  (x^2+1)^{k_0} \, \,du \\  
&=& \sum_{m=0}^M \sum_{j=0}^{\deg p_m} (-1)^m \, p_{m,j}\,\,\int_{-x}^x  \lambda_{j,\,k_0-j,\, \ddot{n}}^{(m)}(u) \,\, f(u) \, du \,\, , 
\end{eqnarray*}
where the convergence $\, \displaystyle \lim_{x\to\infty } Y(x) = \bigl( \, f, \, \tilde{C}\psi_{k_0-j,\, \ddot{n}}\, \bigr)_{(k_0)}$ holds because 
$\tilde{C}\psi _{k_0-s_1,\, \ddot{n}} \in L_{(k_0)}^2(\mathbb{R} )$ and $f \in L_{(k_0)}^2(\mathbb{R} )$.
Then, integrating by parts (which is always applicable to integrals over a finite interval),  
\begin{eqnarray}
Z(x) = W(x) + Y(x) \,\,\,\,\,\,\,\,\, {\rm with} && \,\,\,\,\,\,\,\, 
W(x):=\sum_{m=0}^M\sum_{j=0}^{\deg p_m} p_{m,j}\, w_{m,j}(x) \,\,\,\,\,\,\,\, {\rm and}\\ \nonumber 
w_{m,j}(x):= \sum_{r=0}^{m-1}  (-1)^{m-r-1} \Bigl[ && \left(  \lambda_{j,\, k_0-s_1,\, \ddot{n}}^{(m-r-1)} (x) \right)\cdot \left(({\textstyle\frac{d}{dx}})^r f(x)\right) \\ \nonumber && \hspace{10mm} 
- \left( \lambda_{j,\, k_0-s_1,\, \ddot{n}}^{(m-1-r)} (-x)\right) \cdot \left(({\textstyle\frac{d}{dx}})^r f(-x) \right) \Bigr] \,\, . 
\Label{eqn:integr_parts}\end{eqnarray}
Here, by a recursive use of (\ref{eqn:limit_sm}), 
\begin{eqnarray}
\hspace{1cm}\lim_{x\to\infty } \left(T^{m} Z\right)(x) = \bigl(  \tilde{B}f, \, \psi_{k_0-j,\, \ddot{n}}   \bigr)_{(k_0-j)},  \,\,\, 
\lim_{x\to\infty } \left(T^{m} Y\right)(x) = \bigl( f,  \tilde{C}\psi _{k_0-s_1,\, \ddot{n}} \bigr)_{(k_0)} .
\Label{eqn:limit_adj}\end{eqnarray}

In the following, we will show how the contribution of $W(x)$ in (\ref{eqn:integr_parts}) behaves as $x\to\infty$ under the `blurring' of $x$ by the operator $T$ defined in Section \ref{sec:sm}. From Lemma \ref{lemma:rational-lambda}, there exists a polynomial $R(x)$ of degree not greater than 
$2\nu _{k_0\,\ddot{n}} +m+j+k_0-s_1-r-2$ such that $\displaystyle \, \lambda_{j,\, k_0-s_1,\, \ddot{n}}^{(m-r-1)} (\pm x) = \frac{R(\pm x)}{(x^2+1)^{2\nu _{k_0\,\ddot{n}}+m-r-1}} \,$ where $\nu _{k_0\,\ddot{n}}$ has also been  defined in Lemma 4.2. Hence, with
$Q(x):=x^{2\nu _{k_0\,\ddot{n}} +m+j+k_0-s_1-r-2}\,  R(\frac{1}{x})\, $ which should be a polynomial of $x$ of order not greater than $2\nu _{k_0\,\ddot{n}} +m+j+k_0-s_1-r-2$, 
we have 
$\displaystyle \, \lambda_{j,\, k_0-s_1,\, \ddot{n}}^{(m-r-1)} (\pm x) = (\pm x)^{k_0-m+j-s_1+r}\cdot \, 
\frac{Q(\pm\frac{1}{x})}{\,\, (1+\frac{1}{x^2})^{2\nu _{k_0\,\ddot{n}}+m-r-1} } \,\, $ for $x\ne 0$. Here note that 
\begin{eqnarray*}
\lim_{x\to\pm\infty } \left| \, \frac{Q(\pm\frac{1}{x})}{\,\, (1+\frac{1}{x^2})^{2\nu _{k_0\,\ddot{n}}+m-r-1}}\, \right|  < \infty 
\, , \,\,\,\,\,  
\lim_{x\to\pm\infty } \left| \, x \, \frac{d}{dx} \left(  \frac{Q(\pm\frac{1}{x})}{\,\, (1+\frac{1}{x^2})^{2\nu _{k_0\,\ddot{n}}+m-r-1}}  \,\right)  \right| = \, 0 . 
\end{eqnarray*}
On the other hand, since $f\in D(\tilde{B})\subset L_{(k_0)}^2(\mathbb{R} ) \subset L_{(k_0+j-m-s_1)}^2(\mathbb{R} )$ due to
$k_0+j-m-s_1\le k_0$, Lemma \ref{lemma:conv-sml2} implies that $\displaystyle \lim_{x\to\pm\infty } \bigl(T|\tilde{p\, }|\bigr)(x)=0$ for
 $\tilde{p\, }(x):=x^{k_0+j-m-s_1}f(x)$.   
Then, since $f\in D(\tilde{B})\subset C^M(\mathbb{R} )$, we can apply Lemma \ref{lemma:inv-lhopital} for $\displaystyle \, g(x)=\frac{Q(\frac{1}{x})}{\,\, (1+\frac{1}{x^2})^{2\nu _{k_0\,\ddot{n}}+m-r-1}}\, $ with $\, k_0+j-m-s_1\, $ instead of $k_0$ and with $r$ instead of $n$, where $p\, (x)=\tilde{p\, }(x)$ and $\displaystyle \, h_r(x)=\Bigl(\lambda _{j,\, k_0-s_1,\, \ddot{n}}^{(m-1-r)}(x) \Bigr)\cdot \Bigl(({\textstyle\frac{d}{dx}})^r f(x)\Bigr)$. $\bigl($Here note that $g(x)$ is defined for each fixed $r$, though it depends on $r$.$\bigr)$\,  Its result 
\begin{eqnarray*}
\lim_{x\to\pm\infty } \left( \, T^{r+1} q \, \right)(x) = 0 \,\,\,\,\,\,\,\,\,\,\,\, {\rm for} \,\,\,\,\,\,\,\,\,\,\,\, 
q(x):=\left(\lambda _{j,\, k_0-s_1,\, \ddot{n}}^{(m-1-r)} (x) \right) \cdot \Bigl(({\textstyle\frac{d}{dx}})^r f(x) \Bigr) \,\, 
\end{eqnarray*}
with the definition of $w_{m,j}$ in (\ref{eqn:integr_parts}) implies that $\,\displaystyle \lim_{x\to\pm\infty } \left( \, T^{m} w_{m,j} \right)(x) = 0 \, $ and hence $\,\displaystyle \lim_{x\to\pm\infty } \left( \, T^{m} W \right)(x) = 0 \, $. This convergence, together with the convergences (\ref{eqn:limit_adj}) results in the required statement $\bigl( \, \tilde{B}f, \, \psi_{k_0-j,\, \ddot{n}}  \, \bigr)_{(k_0-j)} = \bigl( f, \, \tilde{C}\psi _{k_0-s_1,\, \ddot{n}} \bigr)_{(k_0)}$, because
$\,\displaystyle \lim_{x\to\pm\infty } \Bigl(\, \left( \, T^{m} Z \right)(x) - \left( \, T^{m} W \right)(x) - \left( \, T^{m} Y \right)(x) \, \Bigr) = 0 \, $ is shown from (\ref{eqn:limit_sm}) and (\ref{eqn:integr_parts}).
\hfill\endproof
\par

\section{Proof of $\bf {\bf \it (iii)}\Longrightarrow 
{\bf \it (iv)}$ under {\bf C1}, {\bf C2}, {\bf C2$^+$}, and {\bf C2.1-C2.3} }
\Label{sec:nenps}
In this section, we will prove that any square-summable vector $\vec{f}$ satisfying $\sum_n b_m^n f_n =0$ corresponds to a true solution in 
${C}^M(\mathbb{R}\backslash {p}_M^{-1}(0)) \cap {\cal H}$ 
of the differential equation  $P(\lambda;x,{\textstyle\frac{d}{dx}}) f=0$, under {\bf C1}, {\bf C2}, {\bf C2$^+$} and {\bf C2.1-C2.3} . %In more precise words, the theorem to be proved is:
In order to show this, we have only to prove the following theorem and the following Theorem. 
\begin{theorem}
\Label{thm:nenps} 
Assume that the quintuplet
$(P(\lambda; x,{\textstyle \frac{d}{dx}}),{\cal H},
\{e_n \, \}_{n=0}^{\infty},
{{\cal H}^\Diamond },
\{e_n^\Diamond \}_{n=0}^{\infty})$
satisfies
Conditions
 {\bf C1}, {\bf C2}, {\bf C2$^+$}, and {\bf C2.1-C2.3}.
Then, any sequence $\vec{f}\in V \cap \ell ^2(\mathbb{Z}^+ )$
satisfies the following.
There exists a function 
$\varphi \in {C}^M(\mathbb{R}\backslash p_M^{-1}(0)) $
such that 
\begin{eqnarray}
P(\lambda;x,{\textstyle\frac{d}{dx}}) \, \varphi (x)=0  
\hbox{ and }
\lim_{N\to\infty } \sum_{n=0}^N f_ne_n(x)=\varphi (x) 
\Label{4-10-1}
\end{eqnarray}
for $\forall x\in \mathbb{R}\backslash p_M^{-1}(0)$.
\end{theorem} 
The proof of this theorem will be constructed in this section. Theorem \ref{thm:nenps} implies that $\displaystyle\sum_{n=0}^N f_ne_n$ converges to a true solution of the ODE as $N\to\infty$ for any $\vec{f}\in V$ in the sense of point-wise convergence 
except at the zero points of  $p_M(x)$. 
Thus, it shows that the statement $(iii)\Longrightarrow (iv)$  holds under the condition in Theorem \ref{thm:equivalence1}.

Especially when $p_M(x)$ has no zero points, Theorem \ref{thm:nenps} guarantees the convergence 
with respect to the ${\cal H}$-norm 
by means of 
the following lemma:
\begin{lemma}
\Label{lemma:conv-norm}
If there exists a function $\varphi \in C^M(\mathbb{R} )$ such that %$\displaystyle P(\lambda;x,{\textstyle\frac{d}{dx}}) \,\varphi (x)=0 $ and 
$\displaystyle \lim_{N\to\infty } \sum_{n=0}^N f_ne_n(x)=\varphi (x) $ holds for any $x\in\mathbb{R} $ for a sequence $\{f_n\}_{n=0}^\infty \in \ell ^2(\mathbb{Z}^+)$, then $\displaystyle \lim _{N\to\infty } \Bigl\|\bigl(\sum_{n=0}^N f_ne_n\bigr)-\varphi \, \Bigr\|_{{\cal H}}=0$.
\end{lemma}
This is just the same as Lemma 3.10 of our preceding paper\cite{paper1}, and the proof is given in that paper.  

% \par\noindent\noindent{\em Proof of Lemma }\ref{lemma:conv-norm}: \quad
% \par \rm 
%Since $\{f_n\}_{n=0}^\infty \in \ell ^2(\mathbb{Z}^+)$ and $\{e_n \, |\, n\in\mathbb{Z}^+\} $ is a CONS of ${\cal H}$, there exists a function $f$ such that  $\displaystyle \lim _{N\to\infty } \Bigl\|\bigl(\sum_{n=0}^N f_ne_n\bigr)-f \Bigr\|_{{\cal H}}=0$. Hence, there exists a subsequence $\{N_\nu \}_{\nu =0}^\infty $ such that $\displaystyle \lim_{\nu \to\infty } \sum_{n=0}^{N_\nu } f_ne_n(x)=f(x) \, $ (a.e.). Therefore, from the trigonometric inequality, $\displaystyle \, |f(x)-\varphi (x)| \le   \lim_{\nu \to\infty } \Bigl( \, \Bigl| \bigl(\sum_{n=0}^{N_\nu } f_ne_n(x)\bigr)-\varphi (x)\, \Bigr| + \Bigl| \bigl(\sum_{n=0}^{N_\nu } f_ne_n(x)\bigr)-f(x)\, \Bigr| \, \Bigr)=0\, $ (a.e.). Therefore, $\|f-\varphi \|_{{\cal H}}=0$,\,  and hence \par\noindent $\displaystyle \lim _{N\to\infty } \Bigl\|\bigl(\sum_{n=0}^N f_ne_n\bigr)-\varphi \,  \Bigr\|_{{\cal H}}\le \lim _{N\to\infty } \Bigl(\, \Bigl\|\bigl(\sum_{n=0}^N f_ne_n\bigr)-f \Bigr\|_{{\cal H}} + \|f-\varphi \|_{{\cal H}} \, \Bigr)=0$.
% \hfill\endproof
% \par \par 
%\par\noindent   Thus, the combination of Theorem \ref{thm:nenps} and Lemma \ref{lemma:conv-norm} shows that the statement $(iii)\Longrightarrow (iv)$  holds under the condition in Theorem \ref{thm:nenps}.

To prove Theorem \ref{thm:nenps}, with the projector $P_{n}$ on $L_{(k_0)}^2(\mathbb{R} )$ to its subspace
${\cal H}^{(n)}:= {\rm span}\bigl(e_0, e_1, \ldots e_n \bigr)$, 
we will analyze the behavior of $\displaystyle \, P_{n} \, y = \sum_{r=0}^n y_r e_r$\, for 
$\vec{y}\in V \cap \ell ^2(\mathbb{Z}^+ )$ 
%$\vec{y}\in V$ 
as $n\to\infty $. Since $\eta =P_{n}f$ is a solution of the inhomogeneous differential equation $P(\lambda;x,{\textstyle \frac{d}{dx})}\eta =g_n$ with $g_n:=P(\lambda;x,{\textstyle \frac{d}{dx})}P_{n}y$ tautologically, we can utilize a kind of 'continuous' correspondence between the inhomogeneous term $g_n$ and the solution $\eta $. There, even though $g_n$ does not converge to $0$ with respect to the $L^2$-norm, the convergence of $\eta$ to a true solution of the homogeneous equation $P(\lambda;x,{\textstyle \frac{d}{dx})}f =0$ can be shown with the help of the characteristic equation  of $N(x,\frac{d}{dx})$ in {\bf C2.2} under some modifications. 

Before giving the proof of Theorem \ref{thm:nenps}, we will provide some preliminaries. First, in order to describe the correspondence between $g_n$ and the $\eta $, we will show some properties of the Green function for the first-order standard form of a $M$th-order differential equation, for any intervals between adjacent zero points of $p_M(x)$, as follows:

 When an inhomogeneous $M$th-order differential equation $\displaystyle\sum_{m=0}^M p_m(x)({\textstyle \frac{d}{dx}})^m \eta =g$ with polynomials $p_m(x)$ ($m=0,1,...,M)$ satisfies the condition that 
$ ^\forall x\in
\tilde{I}, $ $p_M(x)\ne 0$ with an open interval $\tilde{I}=(z,\tilde{z})$  
%^\forall x\in\mathbb{R}, $ $p_M(x)\ne 0$ 
and the  function $g(x)$ is continuous,  we use the following standard form 
\begin{eqnarray}
\frac{d}{dx} \vec{\eta }(x) = {\rm\bf M}(x) \, \vec{\eta } (x) + \, \vec{g}(x)
\Label{eqn:eqgvec}\end{eqnarray}
with the $M$-dimensional vectors  
\begin{eqnarray*}
\displaystyle \left( \vec{\eta } (x)\right)_\ell := \frac{d^\ell }{dx^\ell } \, \eta(x) \,\,\,  (\ell =0,1,...,M-1) \,\, , \,\,\,\,
\left( \vec{g} (x)\right)_\ell  := \left\{ 
\begin{array}{@{\,}ll} 0 & ({\rm if }\,\, 0\le \ell \le M-2) \\ \\ \displaystyle \frac{g(x)}{p_M(x)}  & ({\rm if} \,\, \ell =M-1)
\end{array} \right.
\end{eqnarray*}
and the $M\times M$-matrix 
\begin{eqnarray*}
{\rm\bf M}(x) :=
\left[
\begin{array}{ccccccc}
0 & 1 & 0 & \ldots  & 0 & 0 \\ \\ 
0 & 0 & 1 & \ldots  & 0 & 0 \\ \\ 
  &   &   & \ddots  &   &   \\ \\ 
  &   &   &         &\ddots &   \\ \\ 
0 & 0 & 0 & \ldots  & 0 & 1 \\ \\ 
\displaystyle -\frac{p_0(x)}{p_M(x)}  &\displaystyle -\frac{p_1(x)}{p_M(x)} &\displaystyle -\frac{p_2(x)}{p_M(x)} & \ldots  &\displaystyle -\frac{p_{M-2}(x)}{p_M(x)}&\displaystyle -\frac{p_{M-1}(x)}{p_M(x)}
\end{array}\right] \,\,.
\end{eqnarray*}
Note that $\left[ {\rm\bf M}(x) \right]_{\ell \,  \ell '} = 0\,\,$ 
for $x\in \tilde{I}$
if $\,\,\ell ' \ge \ell +2$. From the existence theorem, the $m$-dimensional vector-valued first-order differential equation (\ref{eqn:eqgvec}) has $M$ linearly independent continuous solutions, because all the elements of ${\rm\bf M}$ are bounded (hence Lipschitz continuity of the right hand side with respect to $\vec{\eta }$ can be derived) and continuous with respect to $x$ and $\vec{g}(x)$ is continuous with respect to $x$ 
for $x\in \tilde{I}$
under the condition that $p_M(x)$ has no real zero.  Therefore, under a choice of the basis vectors, there are $M$ continuous solutions $\vec{\eta }_0(x), \, \vec{\eta }_1(x), \, ... \vec{\eta }_{M-1}(x), \,$ which satisfy the initial conditions $\left( \vec{\eta }_m (\xi) \right)_{\ell } = \delta_{m \, \ell } \,\,\, (\ell =0,1,...,M-1; \, m=0,1,...,M-1)$.  Corresponding to this, consider the following vector-valued standard form of the corresponding homogeneous equation $P(\lambda;x,{\textstyle\frac{d}{dx}}) f = 0$ :
\begin{eqnarray}
\frac{d}{dx} \vec{f} (x)  = {\rm\bf M}(x) \, \vec{f} (x)\, .  
\Label{eqn:eqhomvec}\end{eqnarray}
Here $\vec{f}(x)$ is an $M$-dimensional vector-valued function of $x$ in standard form defined by $\displaystyle\bigl(\vec{f}(x)\bigr)_\ell =({\textstyle\frac{d}{dx}})^\ell f(x)$; it is distinct from $\vec{f}\in \ell ^2(\mathbb{Z}^+ )$ used in other parts of this paper. Let $\vec{f}_0(x), \, \vec{f}_1(x), \, ... \vec{f}_{M-1}(x) \,$ 
($x\in \tilde{I}$)
be its $M$ continuous solutions which satisfy the initial conditions $\left( \vec{f}_m (\xi) \right)_{\ell } = \delta_{m \, \ell } \,\,\, (
\xi \in \tilde{I}; \, 
\ell =0,1,...,M-1; \, m=0,1,...,M-1)$, whose existence is guaranteed in a similar way to the case of (\ref{eqn:eqgvec}).  

Define the $M\times M$-matrix ${\rm\bf \Phi}(x, \xi)$  by $\,\,\displaystyle \left[{\rm\bf \Phi}(x; \xi) \right]_{\ell \, m}:= 
\left( \vec{f}_m (x)\right)_{\ell} \,\,$
for $x, \xi\in \tilde{I}$
, which satisfies $\frac{\partial }{\partial x} {\rm\bf \Phi}(x; \xi) = {\rm\bf M}(x) \, {\rm\bf \Phi}(x; \xi)$ and  ${\rm\bf \Phi}(
\xi; \xi)=I_M$
for $x, \xi\in \tilde{I}$. As is well known, ${\rm\bf \Phi}(x; \xi) $ satisfies the reproducing relation 
\begin{eqnarray}
{\rm\bf \Phi}(x; x') \, {\rm\bf \Phi}(x'; \xi) = {\rm\bf \Phi}(x; \xi)
\hspace{5mm} (x, x', \xi \in \tilde{I}) \, 
\Label{eqn:reproduc}\end{eqnarray}
and another partial differential equation 
\begin{eqnarray}  
\frac{\partial }{\partial \xi} {\rm\bf \Phi}(x; \xi) = - { {\rm\bf \Phi}(x; \xi) \, \rm\bf M}(\xi) \,\, 
\hspace{5mm} (x, \xi\in \tilde{I}).
\Label{eqn:eqadjph}\end{eqnarray}
Partial differentiability of ${\rm\bf \Phi}(x; \xi )$ with respect to $\xi $ is easily shown from the discussion about the difference under an infinitesimal change of $\xi $, because (\ref{eqn:eqadjph}) is derived from the differentiation with respect to $x'$ of both sides of the above reproducing relation (\ref{eqn:reproduc}) and the regularity of the matrices is guaranteed by the linear independence of the columns.

Here, we state a lemma about the higher-order partial derivatives with respect to $\xi$ of  $\left[ \, {\rm\bf \Phi}(x; \xi)\, \right]_{\, 0 \,\, M-1}$, especially at $\xi =x$, which will play an important role later. 
\begin{lemma}
\Label{lemma:differential}
Let $x, \xi\in \tilde{I}$. 
$\left[\, {\rm\bf \Phi}(x; \xi )\, \right]_{\, 0\,\, M-1}$ is partially differentiable with respect to $\xi$ infinitely many times for $\xi \le x$, %and especially for $\ell =0,1,...,M-2$, 
%$\left.\frac{\partial ^m}{\partial \xi ^m} \, \left[ \, {\rm\bf \Phi}(x; \xi ) \, \right]_{\, 0 \, M-1 '} \, \right|_{\xi =x} =\, 0 \,$.
where partial differentiability with respect to $\xi $ for $\xi \le x$ includes the existence of ${\it finite}$ partial differential coefficients from the left at $\xi =x$,
\end{lemma} 
{\em Proof of Lemma }\ref{lemma:differential}: \quad
\rm 
Since $\left[ {\rm\bf M}(x) \right]_{\ell \,  \ell '}\,\,$ is differentiable with respect to $x$ infinitely many times
for $x, \xi\in \tilde{I}$, %and $\,\,\left[ {\rm\bf M}(x) \right]_{\ell \,  \ell '} = 0\,\,$ for $\,\,\ell ' \ge \ell +2$, with (\ref{eqn:eqadjph}) and ${\rm\bf \Phi}(x; x)=I_M$, 
mathematical induction on $m$ by a recursive use of (\ref{eqn:eqadjph}) results in the following (*) for $m\in\mathbb{Z}^+$:
\begin{enumerate}
\item[(*)]
$\frac{\partial ^m}{\partial \xi ^m} \, \left[ \, {\rm\bf \Phi}(x; \xi ) \, \right]_{\, \ell \, \ell '}  \,\,$ are partially differentiable by $\xi$ for $\xi \le x$ 
\hspace{5mm} ($x, \xi\in \tilde{I}$).
%\item[(b)]
%$\left. \frac{\partial ^m}{\partial \xi ^m} \, \left[ \, {\rm\bf \Phi}(x; \xi ) \, \right]_{\, \ell \, \ell '}\, \right|_{\xi =x-0+}\, =0\,\,$ for $\,\, \ell ' \ge \ell +m+1$ \,\, .
\end{enumerate}
\hfill\endproof
With ${\rm\bf \Phi}(x; \xi)$ defined above, as is well known, the relation 
\begin{eqnarray*}
\vec{\eta }_m (x) = {\rm\bf \Phi}(x; \xi) \, \vec{1}_m + \int_{\xi}^x {\rm\bf \Phi}(x; x') \, \vec{g}(x') \, dx' 
\hspace{5mm} (x, \xi\in \tilde{I})
\end{eqnarray*}
holds with $(\vec{1}_m )_{m'} := \delta_{m \, m'}$. Hence, the solution $\vec{\eta }_{\tau }$ of (\ref{eqn:eqgvec}) with the initial conditions $\vec{\eta } (\xi) = \vec{\tau}$ is 
\begin{eqnarray*}
\vec{\eta }_{\vec{\tau }} (x) = {\rm\bf \Phi}(x; \xi) \, \vec{\tau} + \int_{\xi}^x {\rm\bf \Phi}(x; x') \, \vec{g}(x') \, dx' \,\, 
\hspace{5mm} (x, \xi\in \tilde{I}).
\end{eqnarray*}
Hence, if we redefine $g(x)$ by extending its domain to $\mathbb{R}$ by 
\begin{eqnarray*}
g(x)=\left\{\begin{array}{@{\,}ll}  g(x) & (\mbox{if } x\in\tilde{I}) \\ 0 & (\mbox{if } x\in\mathbb{R}\backslash\tilde{I}) ,\end{array} \right. 
\end{eqnarray*}
%In other words, 
under {\bf C2$^+$}, the solution of the inhomogeneous differential equation ${ P(\lambda;x,{\textstyle\frac{d}{dx}})} \eta = g$ 
for $x\in\tilde{I}$ 
with initial conditions $\frac{d^\ell }{dx^\ell } \, \eta (\xi) = \left(\vec{\tau}\right)_\ell \,\,\, (
\xi \in \tilde{I};\,\,  
\ell =0,1,...,M-1)$ can be written in the simple form 
\begin{eqnarray}
\eta_{\vec{\tau}} (x) = \left( \vec{\Phi}(x; \xi) , \, \vec{\tau} \right) + \langle \chi_{\xi ,x} \, , \, g \rangle _{{{\cal H}^\Diamond}}  
\hspace{5mm} (x, \xi\in \tilde{I}).
\Label{eqn:sol_inhom}
\end{eqnarray}
with the vector $\vec{\Phi}(x; x')$ defined by 
\begin{eqnarray*}
\left( \vec{\Phi}(x; x') \right)_\ell := \left[ {\rm\bf \Phi}(x; x') \right]_{\, 0 \,\,\ell } \,\,\, (
x, x' \in\tilde{I};\, 
\ell =0,1,...,M-1) \,\,  
\end{eqnarray*}
and the function 
\begin{eqnarray} 
\chi_{\xi ,u} (x) := 
\frac{ 1_{[\xi, u]}(x) \,\, \left[ \overline{ \, \Phi(u; x) } \, \right]_{\, 0 \,\, M-1} }{
{v}^\Diamond(x)
\,\,\, \overline{p_M(x)} } \,\,\, 
\hspace{5mm} (x, \xi, u\in \tilde{I}),
\Label{eqn:def_chi}
\end{eqnarray}
with ${v}^\Diamond(x)$ in {\bf C2$^+$}, where $1_J(x)$ denotes the indicator function for the interval $J$.

Here, we state a preliminary lemma related to this function,  
where $M$ is the order of $P(\lambda;x,{\textstyle\frac{d}{dx}})\,$ 
and $\tilde{I}=(z, \tilde{z})$ be an open interval in which $p_M(x)$ has no zero points. 
%We have the following:
\begin{lemma}
\Label{lemma:inner-product-unilateral}
Let $\xi \in \tilde{I}$. 
Under {\bf C2$^+$}, {\bf C2.2}, and {\bf C2.3}, for any $u\in \tilde{I}$ greater than $\xi $,
%Under {\bf C2.2}-{\bf C2.4}, for any $u\in\mathbb{R} $ greater than $\xi $, 
$\,\,\, ^\exists K_{\xi ,u}>0$ and $ ^\exists n_c\in \mathbb{Z}^+ \,\,\,$ such that 
$\,\,\, \left| \, \bigl\langle  \, \chi_{\xi ,u} \, , \, e^{\Diamond}_n \bigr\rangle _{{{\cal H}^\Diamond}}\, \right| 
\, \le \displaystyle \,  \frac{K_{\xi ,u}}{n^M} \, \,$ for any $ n \in \mathbb{Z}^+$ greater than  $n_c$. 
\end{lemma} 
\begin{lemma}
\Label{lemma:N-M}
Under {\bf C2.2}, $\bigl(N(x,\frac{d}{dx})\bigr)^M$ can be expressed as the finite sum 
\begin{eqnarray*}
\bigl(N(x,{\textstyle \frac{d}{dx}})\bigr)^M=\sum_{m=0}^M \nu _m(x) ({\textstyle \frac{d}{dx}})^m 
\end{eqnarray*} 
where the functions $\nu _m$ $(m=0,1,\ldots ,M)$ 
%are in 
belong to  
$C^0(\mathbb{R} )$.
\end{lemma}
The proof of this lemma follows easily by mathematical induction on $M$.
\begin{lemma}
\Label{lemma:finite}
Under {\bf C2.2} and {\bf C2.3}, for any real numbers $a$ and $b$ 
in $\tilde{I}$
such that  $a<b$, a function 
$f\in C^M(\tilde{I})$ 
%$f\in C^M(\mathbb{R} )$ 
satisfies the relation  
\begin{eqnarray*}
^\exists C_{a,b}\in \mathbb{R} \mbox{ s.t. } ^\forall n\in\mathbb{Z}^+, \,\, \Bigl| \int_a^b f(x) \,\, \Bigl( \bigl(N(x,\frac{d}{dx})\bigr)^Me^\Diamond_n(x)\Bigr)\,  dx \Bigr| < C_{a,b} .
\end{eqnarray*}
\end{lemma}
{\em Proof of Lemma }\ref{lemma:finite}: \quad
\rm
Under {\bf C2.2}, $\bigl(N(x,{\textstyle \frac{d}{dx}})\bigr)^Me^\Diamond_n$ is well defined. 
From Lemma \ref{lemma:N-M}, with $\mu _m (x):= \nu _m(x)f(x)$
$\in C^M(\tilde{I})\subset C^M[a,b]$, 
\begin{eqnarray*}
&& \int_a^b f(x) \,\, \Bigl( \bigl(N(x,{\textstyle \frac{d}{dx}})\bigr)^Me^\Diamond_n(x)\Bigr)\,  dx = \sum_{m=0}^M \int_a^b f(x) \,\, \Bigl( \nu _m(x) \, (e^\Diamond_n)^{(m)}(x)\Bigr)\,  dx \\  
&=& \sum_{m=0}^M \biggl((-1)^m \int_a^b \mu _m^{(m)}(x) \,\, e^\Diamond_n(x)\,  dx\\
 &&\hspace{10ex}
+\, (-1)^r \sum_{r=0}^{m-1} \bigl( \mu _m^{(r)}(a) \,\, (e^\Diamond_n)^{(m-r-1)}(a) - \mu _m^{(r)}(b) \,\, (e^\Diamond_n)^{(m-r-1)}(b)\bigr) \biggr) \, . 
\end{eqnarray*}
Since $\mu _m$ $(m=0,1,\ldots M-1)$ and $e^\Diamond_n$ $(n\in\mathbb{Z}^+)$ belong to 
$C^M[a,b]$, 
the functions $\mu _m^{(r)}(x)$ and $\tilde{e} _m^{(r)}(x)$ belong to $C^0[a,b]$ for $r\le M$. 
%$C^M(\mathbb{R} )$, 
% he functions $\mu _m^{(r)}$ and $\tilde{e} _m^{(r)}$  are continuous for $r\le M$. 
Hence, under {\bf C2.2} and {\bf C2.3}, all the maxima $\displaystyle M_m^{(r)} := \max_{x\in[a,b]} |\mu _m^{(r)}(x)|$ and $\displaystyle \tilde{A}^{(r)} := \max_{x\in[a,b]} |\tilde{a}^{(r)}(x)|$ $(0\le m\le M$ and $0\le r\le M)$ are finite, 
with $\tilde{a}(x)$ in {\bf C2.3} . From these facts, 
\begin{eqnarray*} 
\Bigl| \int_a^b f(x) \,\, \Bigl( \bigl(N(x,\frac{d}{dx})\bigr)^Me^\Diamond_n(x)\Bigr)\,  dx \Bigr|
 \le \sum_{m=0}^M \Biggl((b-a)M_m^{(m)}\tilde{A}^{(0)} + \sum_{r=0}^{m-1} 2M_m^{(r)}  \tilde{A}^{(m-r-1)}\Biggr)\, ,  
\end{eqnarray*}
where the right hand side is finite and does not depend on $n$.
\hfill\endproof 

{\em Proof of Lemma }\ref{lemma:inner-product-unilateral}: \quad
\rm 
Let $\xi ,u \in \tilde{I}=(z,\tilde{z})$ and let 
$f_{\xi,u}$ be a function in 
$C^M(\tilde{I})$ 
%$C^M(\mathbb{R} )$ 
such that $ f_{\xi,u}(x)= \displaystyle\frac{\left[ \, \overline{\Phi(u; x)}  \, \right]_{\, 0 \,\, M-1}}{\overline{p_M(x)}}  $ for $x\in [\xi,u]$. The existence of $f_{\xi,u}$ is obvious from  the extension of the function to the intervals 
$(z, \xi )$ and $(u, \tilde{z})$ 
%$(-\infty , \xi )$ and $(u,\infty )$ 
by the Taylor expansions of $\displaystyle\frac{\left[ \, \overline{\Phi(u; x)}  \, \right]_{\, 0 \,\, M-1}}{\overline{p_M(x)}}  $ up to the $M$-th order term 
about 
%around 
$x=\xi $ and $x=u$, respectively, because of Lemma \ref{lemma:differential}.

Under {\bf C2$^+$}, {\bf C2.2}, and {\bf C2.3}, since $N(x,\frac{d}{dx})\, e^\Diamond_n(x)=\lambda _n e^\Diamond_n(x)$, 
\begin{eqnarray*}
\lambda _n^M \overline {\bigl\langle  \, \chi_{\xi ,u} \, , \, e^\Diamond_n \bigr\rangle _{{{\cal H}^\Diamond}}}  
= \lambda _n^M  \int_\xi ^u f_{\xi,u}(x) \,e^\Diamond_n(x)\,  dx 
= \int_\xi ^u f(x) \,\, \Bigl( \bigl(N(x,{\textstyle \frac{d}{dx}})\bigr)^Me^\Diamond_n(x)\Bigr)\,  dx\, .
\end{eqnarray*}
This and Lemma \ref{lemma:finite} result in $\displaystyle ^\exists C_{\xi ,u}\in \mathbb{R} $ such that  $|\lambda _n|^M \left|\bigl\langle  \, \chi_{\xi ,u} \, , \, e^\Diamond_n \bigr\rangle _{{{\cal H}^\Diamond}}\right|\le C_{\xi ,u} $. The condition $\displaystyle \liminf_{n\to\infty }\frac{|\lambda _n|}{n}>0$ in {\bf C2.2} implies that there exist an integer $n_c$ and a positive constant $c$ such that $|\lambda _n|\ge c\, n$ may be guaranteed for any $n$ greater than $n_c$. Hence, $\left| \, \bigl\langle  \, \chi_{\xi ,u} \, , \, e^\Diamond_n \bigr\rangle _{{{\cal H}^\Diamond}}\, \right| 
\, \le \displaystyle \,  \frac{C_{\xi ,u}}{(c\, n)^M} \, \,$ for any $ n \in \mathbb{Z}^+$ greater than  $n_c$. With $\displaystyle K_{\xi ,u}:=\frac{C_{\xi ,u}}{c\, ^M} $, the lemma holds.
\hfill\endproof 
\par

Next, as another tool for the proof of the theorem, we will consider the problem of finding the solution of the differential equation $P(\lambda;x,{\textstyle\frac{d}{dx}})\eta =g$ 
for $x\in \tilde{I}=(z,\tilde{z})$ 
under the constraints $\eta (x_j) = t_j \,\,\, (j=0,1,...,M-1)$ for a sequence $x_0 < x_1 < ... <x_{M-1}$
in $\tilde{I}$,  
instead of giving the $M$ initial conditions only at $x=\xi \,\,\, (
z<
\xi <x_0)$. For this problem, define the $M\times M$-matrix ${\rm\bf T}$ by $\left({\rm\bf T}\right)_{j \, m} := f_m (x_j) \,\,\, (j=0,1,...,M-1; \, m=0,1,...,M-1)$ with the solutions $f_m\,\,\,(m=0,1,...,M-1)$ of the homogeneous differential equation $P(\lambda;x,{\textstyle\frac{d}{dx}}) f =0$ 
for $x\in \tilde{I}=(z,\tilde{z})$ 
with the initial conditions $\frac{d^{\,\ell }}{dx^{\,\ell }} \, f(\xi) = \delta_{\,\ell \, m} \,\, (\ell =0,1,...,M-1)$ 
where $z<\xi<x_0$. 
Then the following lemma holds concerning the invertibility of ${\rm\bf T}$:
\begin{lemma}
\Label{lemma:invertibility-T}
When $\displaystyle\, P(x,{\textstyle\frac{d}{dx}})=\sum_{m=0}^M p_m(x) ({\textstyle\frac{d}{dx}})^m $ for polynomials
$p_m$ $(m=0,1,...,M)$ satisfying $\bigl( ^\forall x\in
\tilde{I}
%\mathbb{R}
, \,\, p_M(x)\ne 0\bigr)$, for any $y\in
\tilde{I}
%\mathbb{R} 
$ not smaller than $\xi $, there exists a sequence of finite intervals $[a_0,b_0],\, [a_1,b_1],\, \ldots [a_{M-1},b_{M-1}]$ with $y<a_0<b_0<a_1<b_1<\ldots <a_{M-1}<b_{M-1}<
\tilde{z}
%\infty
$ such that ${\rm\bf T}$ may be invertible when $x_j\in [a_j,b_j]$ $(j=0,1,2,...,M-1)$.
\end{lemma}

{\em Proof of Lemma }\ref{lemma:invertibility-T}: \quad
 \rm 
Define $n\times n$-submatrices ${\rm\bf \tilde{T}}_n(x_0,x_1,...,x_{n-1})$ 
$(n=1,2, ..., M)$\,  by  $\displaystyle\bigl({\rm\bf \tilde{T}}_n(x_0,x_1,...,x_{n-1})\bigr)_{j\, m}:=f_m(x_j)$ 
$(j=0,1,...,n-1;\, m=0,1,...,n-1)$. Then ${\rm\bf \tilde{T}}_M(x_0,x_1,...,x_{M-1})={\rm\bf T}$. Since 
the statement that 
$f_0(x)=0$ for any $x$ 
in $(y,\tilde{z})$ 
%greater than $y$ 
is contradictory to the uniqueness theorem and the initial condition at $x=\xi $, there exists $x_0$ such that 
$y<x_0<\tilde{z}$ 
%$x_0>y$ 
and $f_0(x_0)\ne 0$. Then $\det {\rm\bf \tilde{T}}_1(x_0)=f_0(x_0)\ne 0$. From this initial statement, we can carry out the following mathematical induction:   When  $\det {\rm\bf \tilde{T}}_j(x_0,x_1,...,x_{j-1})\ne 0$, there should exist $x_j$\, 
in $(x_{j-1},\tilde{z})$ 
%  ($(>x_{j-1}$) 
such that  
$\det {\rm\bf \tilde{T}}_{j+1}(x_0,x_1,...,x_j)\ne 0$, because 
$\det {\rm\bf \tilde{T}}_{j+1}(x_0,x_1,...,x_{j-1},x)=0$ for any $x$ 
in $(x_{j-1},\tilde{z})$ 
%greater than $x_{j-1}$ 
would imply $\displaystyle\sum_{m=0}^j c_m(x_0, ...,x_{j-1})\,  f_m(x)=0$ 
for $x>x_{j-1}$ with  $c_j(x_0, ...,x_{j-1})=\det {\rm\bf \tilde{T}}_{j}(x_0,x_1,...,x_{j-1})\ne 0$ which is contradictory to the uniqueness theorem and the initial condition at $x=\xi $. From this mathematical induction, there exists a sequence $y<x_0 < x_1 < ... <x_{M-1}
<\tilde{z}
$ such that $\det {\rm\bf \tilde{T}}_M(x_0,x_1,...,x_{M-1})\ne 0$ \, i.e. $\det {\rm\bf T} \ne 0$. 

Next, from the conditions for $P(x,{\textstyle\frac{d}{dx}})$ and the existence theorem,  
$\det {\rm\bf T} = \det {\rm\bf \tilde{T}}_M(x_0,x_1,...,x_{M-1})$ is $M$-times continuously partially  differentiable with respect to $x_j$ ($j=0,1,...,M-1$) 
in $\tilde{I}$ 
and moreover totally differentiable 
$\tilde{I}^M$, and hence it is locally Lipschitz continuous 
there. 
Therefore, with the conventional vector notation $\vec{x}\in \mathbb{R} ^M$ defined by $(\vec{x})_j=x_j$  ($j=0,1,...,M-1$), if $\det {\rm\bf \tilde{T}}_n(x_0,x_1,...,x_{n-2},x)\ne 0$, there exists a neighborhood $\displaystyle U_\epsilon (\vec{x})=\bigl\{ \vec{u}\,\bigl|\, \|\vec{u}-\vec{x}\|<\epsilon \bigr\}$\,  $(\epsilon >0)$ in 
$\tilde{I}^M$ 
%$\mathbb{R} ^M$ 
such that $\det {\rm\bf \tilde{T}}_n(u_0,u_1,...,u_{n-1})\ne 0$ for any $\vec{u}\in U_\epsilon (\vec{x})$. 
Since $\{\vec{u} \, |\, u_j\in [x_j-\delta _j,\, x_j+\delta _j] \,\, (j=0,1,...,M-1)\}\subset U_\epsilon (\vec{x})$ holds at least for $0< \delta _j< \frac{\epsilon}{\sqrt{M}}$\, $(j=0,1,...,M-1)$, the lemma holds with $a_j:=x_j-\delta _j$ and $b_j:=x_j+\delta _j$ \,  (where 
$z<b_j<a_{j+1}<\tilde{z}$
%$b_j<a_{j+1}$ 
is satisfied for an appropriate choice of sufficiently small $\delta _j$ and $\delta _{j-1}$).
\hfill\endproof

Under the existence of a sequence with invertible ${\rm\bf T}$ guaranteed by this lemma, we have another lemma with the definition of the vector $\vec{b}_g$ defined by 
\begin{eqnarray} 
(\vec{b}_g)_j := \langle \chi_{\xi ,x_j} \, , \, g \rangle _{{{\cal H}^\Diamond}}\,\,
\,\,\,\,\,\,\, (j=0,1,...,M-1) .
\Label{eqn:def_b}\end{eqnarray}
\begin{lemma}
\Label{lemma:sol-constr}
When the sequence 
$z<x_0 < x_1 < ... <x_{M-1}<\tilde{z}$ 
%$x_0 < x_1 < ... <x_{M-1}$ 
is chosen so that ${\rm\bf T}$ is invertible, the solution of the inhomogeneous differential equation $P(\lambda;x,{\textstyle\frac{d}{dx}})\eta =g$ 
for $x\in\tilde{I}$ 
under the constraints $\eta (x_j) = t_j \,\,\, (j=0,1,...,M-1)$ \par\noindent  $($where 
$z<\xi < x_0 < x_1 < .... <x_{M-1}<\tilde{z})$
%$\xi < x_0 < x_1 < .... <x_{M-1})$
is 
\begin{eqnarray*}
\eta_{{\rm\bf T}^{-1} ( {\vec{t}} - \vec{b}_{g} )} (x) = \left( \vec{\Phi}(x; \xi) , \, {\rm\bf T}^{-1} ( {\vec{t}} - \vec{b}_g ) \, \right) + \langle \chi_{\xi ,x} \, , \, g \rangle_{{{\cal H}^\Diamond}}
\end{eqnarray*}
with the vector ${\vec{t}}\in \mathbb{R} ^M$ defined by $\left({\vec{t}}\right)_j =t_j\,$ \, $(j=0,1,...,M-1)$.
\end{lemma} 
{\em Proof of Lemma }\ref{lemma:sol-constr}: \quad
Since the homogeneous differential equation $P(\lambda;x,{\textstyle\frac{d}{dx}}) f =0$ 
for $x\in\tilde{I}$ 
is a special case of the inhomogeneous differential equation with $g=0$, from (\ref{eqn:sol_inhom}), the solution of the homogeneous differential equation $P(\lambda;x,{\textstyle\frac{d}{dx}}) f =0$ 
for $x\in\tilde{I}$ 
with the initial conditions $ \frac{d^\ell }{dx^\ell } \, f(\xi) = \left( \vec{\tau} \right)_\ell \,\,\, (\ell =0,1,...,M-1)$ is 
$f_{\vec{\tau}} (x_j) = \left( \vec{\Phi}(x_j; \xi) , \, \vec{\tau} \right)\,$. 
As special cases, we have $f_{m} (x_j) = \left( \vec{\Phi}(x_j; \xi) , \, \vec{1}_{m} \right) $ with the vector $\vec{1}_m$ defined by $\left(\vec{1}_m\right)_\ell :=\delta _{m \, \ell }\,$. Define the $M$-dimensional vector $\vec{t}_{\vec{\tau}}$ such that  $\left(\vec{t}_{\vec{\tau}}\right)_j=\eta_{\vec{\tau}}(x_j)   (j=0,1,...,M-1)$. Since $\,\eta_{\vec{\tau}} (x_j)-\left( \vec{b}_{g} \right)_j=f_{\vec{\tau}} (x_j)$, from the above relations, we have 
\begin{eqnarray*}
\vec{t}_{\vec{\tau}} - \vec{b}_{g} 
&=& \left[\begin{array}{@{\,}ll} 
\left( \vec{\Phi}(x_0; \xi) , \, \vec{\tau} \right) \\ 
\left( \vec{\Phi}(x_1; \xi) , \, \vec{\tau} \right) \\ 
\,\,\,\,\,\,\,\,\,\,\,\,\,\,\, . \\ 
\,\,\,\,\,\,\,\,\,\,\,\,\,\,\, . \\ 
\left( \vec{\Phi}(x_{M-1}; \xi) , \, \vec{\tau} \right) 
\end{array}\right]
= \sum_{m=0}^{M-1} \left(\vec{\tau}\right)_m \, 
\left[\begin{array}{@{\,}ll} 
\left( \vec{\Phi}(x_0; \xi) , \, \vec{1}_m \right) \\ 
\left( \vec{\Phi}(x_1; \xi) , \, \vec{1}_m \right) \\ 
\,\,\,\,\,\,\,\,\,\,\,\,\,\,\, . \\ 
\,\,\,\,\,\,\,\,\,\,\,\,\,\,\, . \\ 
\left( \vec{\Phi}(x_{M-1}; \xi) , \, \vec{1}_m \right) 
\end{array}\right]\\ 
&=& \sum_{m=0}^{M-1} \left(\vec{\tau}\right)_m \, 
\left[\begin{array}{@{\,}ll} 
f_m(x_0) \\ 
f_m(x_1) \\ 
\,\,\,\,\,\,\,\,\,\,\,\,\,\,\, . \\ 
\,\,\,\,\,\,\,\,\,\,\,\,\,\,\, . \\ 
f_m(x_{M-1}) \\ 
\end{array}\right]
= {\rm\bf T} \, \vec{\tau} . 
\end{eqnarray*}
Hence, we can show that 
the function 
\begin{eqnarray*}
\eta_{{\rm\bf T}^{-1} ( {\vec{t}} - \vec{b}_{g_n} )} (x) = \left( \vec{\Phi}(x; \xi) , \, {\rm\bf T}^{-1} ( {\vec{t}} - \vec{b}_{g_n} ) \, \right) + \langle \chi_{\xi ,x} \, , \, g \rangle_{{{\cal H}^\Diamond}}
\end{eqnarray*}
is the solution of $P(\lambda;x,{\textstyle\frac{d}{dx}}) \eta = g $
for $x\in\tilde{I}$ 
satisfying the constraints $\eta (x_j) = ({\vec{t}})_j $ \par\noindent $ (j=0,1,...,M-1)$ for the sequence 
$z<\xi < x_0 < x_1 < .... < x_{M-1} < x_M <\tilde{z}$, 
%$\xi < x_0 < x_1 < .... < x_{M-1} < x_M $, 
where the uniqueness of the solution satisfying these constraints has been shown  also.
\hfill\endproof

By using these preliminaries, now we are able to construct the proof of Theorem \ref{thm:nenps} as follows; 

{\em Proof of Theorem }\ref{thm:nenps}: \quad

Suppose that $\vec{y} \in V \cap \ell ^2(\mathbb{Z}^+ )$
has no function $\varphi \in {C}^M(\mathbb{R}\backslash p_M^{-1}(0)) $ 
satisfying (\ref{4-10-1}).
Then, 
the basis syetem $\{e_n\, |\, n\in\mathbb{Z}^+ \} $ is a CONS of ${\cal H}=L_{(k_0)}^2(\mathbb{R} )$, the function $\displaystyle y:=\sum_{n=0}^\infty y_ne_n$ belongs to ${\cal H}$. With the projector $P_{n} \, $ on  ${\cal H}$ to the subspace  $\,{\cal H}^{(n)}:= {\rm span}(e_0, e_1, ... , e_n)$, the convergence 
$\displaystyle \lim_{n\to\infty }\|P_{n} \, y - y\|_{({k_0^\Diamond})}=0$ holds.   
Hence, there exists a subsequence $\{n_\nu \}_{\nu =0}^\infty $ such that  $\displaystyle\lim_{\nu\to\infty }\bigl(P_{n_\nu } \, y\bigr)(x)=y(x)$ (a.e.). 
 
Therefore, from the assumption that $\vec{y} \in (V^{\ell ^2}\backslash V) \subset V$ and Lemma \ref{lemma:invertibility-T}, without loss of generality, we can show the existence of a sequence 
$\xi < x_0 < x_1 < .... < x_{M-1} < x_M $ $\bigl($where $M$ is the order of the differential equation $P(\lambda;x,{\textstyle\frac{d}{dx}})f=0\bigr)$ %and an  initial point $\xi \,\, (<x_0 -1)$
 satisfying the following conditions (a)-
(d):
%(c):
\begin{enumerate}
\def\labelenumi{(\alph{enumi})}
\item 
For $j=0,1,...,M$, the limits $\displaystyle \lim_{\nu\to\infty } \bigl(P_{n_\nu } \, y\bigr )(x_j)$ exist and 
$\displaystyle \lim_{\nu\to\infty } \bigl(P_{n_\nu } \, y\bigr )(x_j) = t_j$. 
\item
The $M\times M$-matrix ${\rm\bf T}$ is invertible under the definition by $({\rm\bf T})_{j \, m} = f_m(x_j)$ \,  $(j=0,1,...,M-1;\, $  $m=0,1,...,M-1)$ for the continuous solutions $f_m $ of $P(\lambda;x,{\textstyle\frac{d}{dx}}) f = 0 $ with the initial conditions  $(\frac{d}{dx})^\ell \, f(\xi) = \delta_{\ell \, m}$ $ (\ell =0,1,...,M-1)$. 
\item
$f(x_M)\ne t_M$ for the true continuous solution $f(x)$ of $P(\lambda;x,{\textstyle\frac{d}{dx}})f=0$ which satisfies $f(x_j)=t_j \,\,\, (j=0,1,...,M-1)$.
\item 
$p_M(x)$ has no zero point in $[\xi ,x_M]$.
\end{enumerate}

Define $\displaystyle g_n (x) := \left( P(\lambda;x,{\textstyle\frac{d}{dx}}) ( P_{n} \, y ) \right)(x)\, $. Then $g_n (x)$ belongs to $C^\infty (\mathbb{R})$  because {\bf C2} implies that it can be expressed as a finite sum of the basis functions $e^\Diamond_n(x)$ of $\tilde{{\cal H}}$. 
Moreover, define the $M$-dimensional vectors $\vec{t}$ and $\vec{t}^{\,\, \nu }$ by $\left(\vec{t}\right)_j:=t_j$ 
and 
$\displaystyle \left(\vec{t}^{\,\, \nu }\right)_j:=\bigl(P_{n_\nu } \, y\bigr)(x_j)$ $ (j=0,1,...,M-1)$, respectively. 
Since the condition (d) requires that $\xi, x_0, x_1, \ldots , x_M$ belong to an open interval $\tilde{I}$ in which $p_M(x)$ has no zero points, from the 
definition, $P_{n_\nu } \, y (x)$ is just the solution of the inhomogeneous differential equation  $P(\lambda;x,{\textstyle\frac{d}{dx}}) \eta(x) = g_{n_\nu }(x)$ for $x\in \tilde{I}$ 
under the constraints 
$\eta (x_j)=\left(\vec{t}^{\,\, \nu }\right)_j \,\, (j=0,1,...,M-1)\,$.    
Therefore, from Lemma \ref{lemma:sol-constr}, 
\begin{eqnarray*}
\bigl(P_{n_\nu } \, y\bigr) (x_M) =  \left( \vec{\Phi}(x_M; \xi) , \, {\rm\bf T}^{-1} ( \vec{t}^{\,\, \nu } - 
\vec{b}_{g_{n_\nu }} 
%\vec{b}_{g_{\xi ,n_\nu }} 
) \, \right) + 
\langle \chi_{\xi ,\, x_{_{M}}} \, , \, 
g_{n_\nu } 
%g_{\xi;n_\nu } 
\rangle _{{{\cal H}^\Diamond}}
\end{eqnarray*}
where the function $\chi_{\xi ,\, x_M}$ and the vector $\vec{b}_g$ have been defined in (\ref{eqn:def_chi}) and (\ref{eqn:def_b}), respectively. 
On the other hand, with $g=0$ in the same lemma, similarly we have 
\begin{eqnarray*}
f(x_M) = \left( \vec{\Phi}(x_M; \xi) , \, {\rm\bf T}^{-1} \, {\vec{t}}  \,\,\, \right).
\end{eqnarray*}
Hence, 
\begin{eqnarray*}
\bigl( P_{n_\nu } \, y \bigr)(x_M)-\bigl( f \bigr) (x_M) = \left( \, \vec{\Phi}(x_M; \xi) , \, {\rm\bf T}^{-1} \bigl( \, (\vec{t}^{\,\, \nu } -\vec{t}\, ) - 
\vec{b}_{g_{n_\nu }}
%\vec{b}_{g_{\xi ,n_\nu }} 
\bigr) \, \right) + \langle \chi_{\xi ,\, x_{_{M}}} \, , \, 
g_{n_\nu }
%g_{\xi ,n_\nu }
\rangle _{{{\cal H}^\Diamond}}  .
\end{eqnarray*}
From the definitions, the limits $\,\,
\vec{b}_{g_{n_\nu }} 
%\vec{b}_{g_{\xi ,n_\nu }} 
\to 0\,\,$ and $\,\,\langle  \chi_{\xi ,\, x_{_{M}}} \, , \, 
g_{n_\nu }
%g_{\xi;n_\nu }
\rangle\to 0\,\, $ as $\,\, \nu \to\infty \,$  
holds if the convergence $\displaystyle \lim_{n\to\infty }\langle  \chi_{\xi ,\, x_j}\, , \, 
g_n
%g_{({k_0^\Diamond}) n }
\rangle=0\,\,$ 
holds for $\,\, j=0,1,...,M\,\,$.

Now, we will prove this convergence,  as follows:
From {\bf C2}, when $n\ge 2\ell $, it is easily shown that  
$\displaystyle \,  g_n(x) =  \sum_{r =n-\ell +1}^{n+\ell } \, \langle g_n , \, e^\Diamond_r  \rangle _{{{\cal H}^\Diamond}}\,\, e^\Diamond_r (x) \, $, \, 
because $\displaystyle \sum _{r=n-\ell }^{\ell } b_m^r \, y_r =0$\, $(m\in\mathbb{Z}^+)$ holds for 
$\vec{y}\in V \cap \ell ^2(\mathbb{Z}^+ )$
and hence $\displaystyle \langle g_n, \, e^\Diamond_m \rangle _{{{\cal H}^\Diamond}}= \sum_{r=m-\ell }^{\max (n, \, m+\ell )} b_m^r\, y_r$ vanishes when $m+\ell \le n$. 
Hence, from Lemma \ref{lemma:inner-product-unilateral}, when $n\ge \ell +1$, 
%\begin{eqnarray}\begin{array}{@{\,}ll}x{cc}\displaystyle 
\begin{eqnarray*}
\left| \, \langle \chi_{\xi ,\, x_j} \, , \, g_n \rangle _{{{\cal H}^\Diamond}}\, \right| 
 = \left|   \sum_{r =n-\ell +1}^{n+\ell } \, \langle \chi_{\xi ,\, x_j} \, , \, e^\Diamond_r  \rangle _{{{\cal H}^\Diamond}}\,\, \overline{\langle g_n , \, e^\Diamond_r  \rangle _{{{\cal H}^\Diamond}}} \, \right| 
%\\ \\ \displaystyle
 \le    K_{\xi ,x_j} \!\!\!\!\! \sum_{r =n-\ell +1}^{n+\ell }
 \frac{\left| \, \langle g_n , \, e^\Diamond_r  \rangle _{{{\cal H}^\Diamond}}  \right|}{r^M}  .
%\end{array}
\Label{eqn:ineqps}\end{eqnarray*}
Here, for $n-\ell +1\le r \le n+\ell $,
\begin{eqnarray}
\frac{\left| \, \langle g_n , \, e^\Diamond_r  \rangle _{{{\cal H}^\Diamond}} \, \right|}{r^M}
=   \frac{1}{r^M} \left| \,  \sum_{\ell =-\ell }^{n-r} 
 b_r ^{r +\ell} \, y_{r +\ell } \,  \, \right| 
\, \le \frac{1}{r^M}  \sum_{\ell =-\ell }^{n-r}   
\left| \, \, b_r ^{r+\ell } \, \right| \,\, \cdot \,\,  
\left| \,y_{r+\ell } \, \right|.
\Label{eqn:ineqvec1}\end{eqnarray}
From {\bf C2.1}, the finite supremum $\displaystyle\,\, K'  \, := \, \sup_{r\in \mathbb{Z}^+\backslash \{0\}, \,\, n\in\mathbb{Z}^+} \frac{| \, b_r ^{n} |}{r^M} \, $ exists.\,  Hence, from (\ref{eqn:ineqvec1}), for   $n-\ell +1\le r \le n+\ell $,  we have 
\begin{eqnarray}
\frac{\left| \, \langle g_n , \, e^\Diamond_r  \rangle _{{{\cal H}^\Diamond}} \, \right|}{r^M}
\, \le \,  K' \, \sum_{\ell =-\ell }^{n-r} \, \left| y_{r+\ell } \, \right| 
 \le \, K' \sqrt{\,  2\,\ell \, \sum_{\ell =-\ell }^{n-r} \, 
\left| y_{r+\ell } \, \right|^2 \, }
\Label{eqn:ineqvec2}\end{eqnarray}
where the last inequality is derived from the Schwartz inequality. 
From the inequalities (\ref{eqn:ineqps}) and  (\ref{eqn:ineqvec2}), for $n\ge 2\,\ell $, we have the inequality 
\begin{eqnarray*}
\left| \, \langle \chi_{\xi ,\, x_j} \, , \, g_n \rangle _{{{\cal H}^\Diamond}}\, \right| 
\,  \le \, 2\,\ell \,  K_{\xi ,x_j}  \, K'\, \, \sqrt{ 2\,\ell  \sum_{r =n-2\ell +1}^\infty \, 
\left| y_{r+\ell } \, \right|^2 \, }\, .
\end{eqnarray*}
Since  $\vec{y}$ is square-summable, $\displaystyle\sum_{r =n-2\ell +1}^\infty \, 
\left| y_{r+\ell } \, \right|^2 \, \to 0 \,\,$ as $\,\, n\to\infty$. Therefore, 
\begin{eqnarray}
\lim_{n\to\infty} \, \langle \chi_{\xi ,\, x_j} \, , \, g_n \rangle _{{{\cal H}^\Diamond}}\, = 0 \,\,\,\,\,\,\,\, 
\mbox{ and hence } \,\,\,\,\,\,\,\, 
\lim_{\nu \to\infty} \, \langle \chi_{\xi ,\, x_j} \, , \, g_{n_\nu } \rangle _{{{\cal H}^\Diamond}}\, = 0  .
\Label{eqn:ineqpsc}\end{eqnarray}
Thus, we have proved that $\displaystyle\lim_{\nu \to\infty} \langle \chi_{\xi ,\, x_j} \, , \, g_{n_\nu } \rangle _{{{\cal H}^\Diamond}}= 0\,\,$, i.e.   $\,\,\displaystyle\lim_{\nu n\to\infty } 
\vec{b}_{g_{n_\nu }}
%\vec{b}_{g_{\xi;n_\nu }}
= 0\,\,$ and $\displaystyle\,\,\lim_{\nu \to\infty} 
\langle \chi_{\xi, x_M} \, , \, g_{n_\nu } \rangle _{{{\cal H}^\Diamond}}
%\langle \chi_{({k_0^\Diamond})\, M} \, , \, g_{\xi;n_\nu } \rangle _{{{\cal H}^\Diamond}}
=0$. 
These convergences, together with the convergence
$\displaystyle \lim_{\nu \to\infty} \vec{t}^{\,\, \nu } = \vec{t}$ which is identical to (a), 
lead us to the conclusion that  $\displaystyle \lim_{\nu \to\infty} \bigl(P_{n_\nu } \, y \bigr)(x_M) = f(x_M)$, 
which is contradictory to (c).  
Therefore, the assumption that
$\vec{y} \in V \cap \ell ^2(\mathbb{Z}^+ )$
does not satisfies the conclusion of 
Theorem \ref{thm:nenps}.
That is, we obtain Theorem \ref{thm:nenps}.

\hfill\endproof
\par

\section{Discussion}

The $m$-th order derivatives ($m=1,2,\ldots ,M-1$) of the eigenfunction $f$ in $L_{(k_0)}^2(\mathbb{R} )$ do not always belong to $L_{(k_0)}^2(\mathbb{R} )$. For example, for the differential operator,   
\begin{eqnarray*} 
{P}(x,{\textstyle\frac{d}{dx}})= (3x^2+1)^2\bigl(\frac{d}{dx}\bigr)^2  + 6x(3x^2+1)\, \bigl({\frac{d}{dx}}\bigr)   - (3x^2+1)^4-18x^2   \, , 
\end{eqnarray*} 
an eigenfunction associated with the eigenvalue $-6$ is $\displaystyle\,f(x)=\frac{1}{3x^2+1}\, \cos (x^3+x)\,$, which belongs to 
$ L^2(\mathbb{R} )$. However, $\displaystyle\frac{d}{dx}f(x) = -\sin (x^3+x)-\frac{6x}{(3x^2+1)^2 }\cos (x^3+x) \notin L^2(\mathbb{R} )$. There are many similar examples. In order to discuss regularity in the framework based on the Sobolev space, some transformation or some change of variable is necessary, for these cases.

Our proof does not require any assumption about whether or not the $m$-th order derivatives ($m=1,2,\ldots ,M-1$) of the eigenfunction $f$ belong to $L_{(k_0)}^2(\mathbb{R} )$. Hence, it can be used to show regularity for these cases without any transformation or modification.

\section{Conclusion}

When an operator 
$\widetilde{A}_{R,L_{(k_0)}^2(\mathbb{R} )}$
%$\tilde{A}_P$
is defined on $C^M(\mathbb{R} )\cap L_{(k_0)}^2(\mathbb{R} )$ 
as the action of an $M$-th order differential operator 
$R(x,\frac{d}{dx})$ 
%$P(x,\frac{d}{dx})$ 
with rational coefficient functions,
we have proved the regularity of the eigenfunctions of the 
closed extension on $L_{(k_0)}^2(\mathbb{R} )$ of 
the operator 
$\widetilde{A}_{R,L_{(k_0)}^2(\mathbb{R} )}$ except at the singular points of the ODE  $R(x,\frac{d}{dx})f=\lambda f$,  
%$\tilde{A}_P$
% except at the zero points of the coefficient function of the highest order, 
This derivation does not require any assumptions about the  $m$-th order derivatives ($m=1,2,\ldots ,M-1$) of the eigenfunction. In particular, for $k_0=0$, we have proved it for the usual $L^2(\mathbb{R} )$. 

The proof has been constructed in two steps: 
The first step has proved regularity in a general framework under several assumptions. The second step has shown that 
the above mentioned operator satisfies these required assumptions.

In the first step, 
the differential operator is treated as an operator from 
a dense subset of a Hilbert space ${\cal H}$ to another Hilbert space ${\cal H}^\Diamond $ which contains ${\cal H}$ (in the sense of sets), 
and this operator can be represented in matrix form with respect to appropriate basis systems of 
${\cal H}$ and ${\cal H}^\Diamond $.
The proof in this framework has been based on 
the implications ${\it(i)\Longrightarrow (ii)\Longrightarrow (iii)\Longrightarrow (iv)}$,
in Theorem \ref{thm:equivalence1} with a more general framework:  
${\it (i)}$ the kernel of the closed extension of the operator on ${\cal H}$ defined as the action of the differential operator, 
${\it (ii)}$ the kernel of the closed extension of the operator from a dense subspace of ${\cal H}$ to ${\cal H}^\diamond $ defined as the action of the differential operator,  
${\it (iii)}$ the space of square-summable number sequences satisfying the simultaneous linear equations corresponding to the matrix representation of one of the above two operators  and 
${\it (iv)}$ the space of `regular' solutions of the differential equation 
which are continuously differentiable $M$ times at any points except the zero points of the coefficient function of the highest order. 
This general framework was used also for an integer-type algorithm for solving higher order homogeneous linear ordinary equations in our preceding paper~\cite{paper1}. 

In the second step, we have shown that the choices ${\cal H}^\Diamond =L_{(k_0^\Diamond )}^2(\mathbb{R} )$ and the basis function systems in (\ref{eqn:choice-bs}) satisfy the  conditions required for the framework in the first step.

The proofs in the two steps have easily been constructed except for two points; 
one is the proof of $\it(iii)\Longrightarrow (iv)$ in the first step and 
the other is the proof of the fact that the choices satisfy condition {\bf C3} in  the second step. 
For the latter point, we have developed a kind of smoothing operator, as a tool. 
Our proof of ${\it (iii)\Longrightarrow (iv)}$ 
has been made by means of a modified kind of continuity of the solutions of an inhomogeneous equation with respect to the inhomogeneous term.

Similar proofs of regularity may be possible even for other choices of function spaces and basis systems satisfying the conditions in this paper or similar type of conditions, which will be a topic for future research. 

\section*{Acknowledgments}
MH was partially supported by MEXT through a Grant-in-Aid for Scientific Research on the Priority Area "Deepening and Expansion of Statistical Mechanical Informatics (DEX-SMI)", No. 18079014 and
a MEXT Grant-in-Aid for Young Scientists (A) No. 20686026. 
The Center for Quantum Technologies is funded by the Singapore Ministry of Education
and the National Research Foundation as part of the Research Centres of Excellence
programme. 
The authors are grateful to
Prof. Hideo Kubo, Tohoku University, for his helpful advices and comments.  

\appendix

\section{Relationship with Fourier series}
\Label{app:fou}

The basis systems used in our methods are closely related to Fourier series expansions of functions defined on the interval $[-\pi,\pi ]$, by the change of variable $\theta = 2 \arctan x \,\, (\, {\rm or} \,\,\, x = \tan \frac{\theta}{2}\, )$. By this change, 
\begin{eqnarray}\begin{array}{@{\,}ll}\displaystyle 
dx = \frac{1}{2} \sec^2 \frac{\theta}{2} \,\,\, d\theta , \,\,\,\,\,\, d\theta = \frac{2}{x^2+1} \,\, dx \, ,  
\,\,\,\,\, 
x^2 + 1 = \sec^2 \frac{\theta}{2} \, ,  
\,\,\,\,\, 
\arg (x\pm i) = \mp \frac{1}{2} \, ( \theta - \pi ) \, ,
\\ \\ \displaystyle 
x\pm i = \pm i \,\,\, e^{\mp i \frac{\theta}{2}} \,\, \sec \frac{\theta}{2} \,\, , 
\,\,\,\,\, 
\frac{x-i}{x+i} = - \,\,e^{i\theta} \, , 
\end{array}\Label{eqn:properties_theta}\end{eqnarray}
and then 
\begin{eqnarray*}
\psi_{k_0,\ddot{n}}(\tan \frac{\theta}{2}) = i^{\, k+1} (-1)^{\ddot{n}} \,\,\, e^{i(\ddot{n}+\frac{k+1}{2})\theta } \,\, \cos^{k+1} \frac{\theta}{2} \,\, . 
\end{eqnarray*}
Here, define 
\begin{eqnarray}
\,\,\,\,\,\,\,\,\,\,\, \tilde{f}(\theta ) := \left\{\begin{array}{@{\,}ll} \displaystyle
\frac{1}{\sqrt{2}} \,\, e^{-\frac{i(k+1)}{2}(\theta +\pi)} \,\, \left| \, \sec^{k+1} \frac{\theta}{2} \, \right| \,\,\, f(\tan \frac{\theta}{2}) &(\mbox{ if } -\pi < \theta < \pi) \\ \\ 0 &  (\mbox{ if } \theta = \pm\pi)\,\, . \end{array}\right.
\Label{eqn:def_tildef}\end{eqnarray}
Then, from the relations (\ref{eqn:properties_theta}), we have a kind of isometric relation 
\begin{eqnarray}
(f, \, g)_{(k_0)} =  \int_{-\pi }^{\pi }  \tilde{f}_k (\theta ) \,\, \overline{\tilde{g}_k (\theta )} \,\, d\theta  \,\, . 
\Label{eqn:ip_redefinedf}\end{eqnarray}
% \begin{rem}\begin{indention}{.8cm}\rm
% Note that $x\to\pm\infty$ is corresponding to $\theta\to\pm\pi$. Hence singularities at $\theta =\pm\pi$ are allowed if we do not require any boundary condition at $x=\pm\infty$. To avoid them, we should define the inner product by the integration in the interval excluding the point $\theta =\pm\pi$. This is an important difference from the case of Fourier series where the inner product is defined by the integration in $[-\pi ,\, \pi]$. Singularities at $\theta =\pm\pi$ are often corresponding to (non-square-summable) pseudo-solutions~\cite{paper1}  of the simultaneous linear equations, which can be removed by our method, as is shown in numerical examples ~\cite{paper1} ~\cite{paper3}.
%\\end{indention}\end{rem}
The relation (\ref{eqn:properties_theta}) and the definitions (\ref{eqn:def-psi}) (\ref{eqn:def_tildef}) result in 
\begin{eqnarray}
\tilde{\psi}_{k,\ddot{n}} (\theta ) = \left\{\begin{array}{@{\,}ll}\displaystyle  \frac{(-1)^{\ddot{n}}}{\sqrt{2}} \, \, \,e^{i\ddot{n}\theta } & (\mbox{ if } -\pi < x < \pi) \\ \\ \displaystyle 
0 & (\mbox{ if } x=\pm\pi ) \end{array}\right.  , 
\Label{eqn:psi_exp}\end{eqnarray}
where we have the characteristic equation  
\begin{eqnarray}
\frac{d}{d\theta} \, \tilde{\psi}_{k,\ddot{n}}(\theta ) =  i\ddot{n}\, \tilde{\psi}_{k,\ddot{n}}(\theta )
\,\,\,\,\,\,\,\,\,\,\,\, (-\pi < \theta < \pi ) 
\Label{eqn:char-eq_exp}\end{eqnarray}
which corresponds to the characteristic equation (\ref{eqn:char-eq}) in the $x$-coordinate in Section \ref{sec:pp}. From (\ref{eqn:psi_exp}), the expansion of $f\in L_{(k_0)}(\mathbb{R} )\,$ with respect to the biorthonormal basis system $\{\sqrt{\frac{1}{\pi }}\, \psi _{k,\ddot{n}} \, | \, \ddot{n}\in\mathbb{Z}\}$, 
\begin{eqnarray}
\,\,\,\,\,\,\,\,\,\,\, f(x)=\frac{1}{\sqrt{\pi}} \sum_{\ddot{n}=-\infty }^\infty \ddot{f}_{\ddot{n}} \,\psi _{k,\ddot{n}}(x) 
\,\,\,\,\,\,\,\,\,\, {\rm with } \,\,\,\,\,\,\,\,\,\, 
\ddot{f}_{\ddot{n}} =\frac{1}{\sqrt{\pi}} \int_{-\infty }^\infty f(x) \,\, \overline{ \psi _{k,\ddot{n}}}\, (x^2+1)^k dx
\Label{eqn:exp_psi}\end{eqnarray}
just corresponds to %the Fourier series expansion 
\begin{eqnarray}
\tilde{f}(\theta ) = \sum_{\ddot{n}=-\infty }^\infty \tilde{F}_{\ddot{n}} \, e^{i\ddot{n}\theta }
\,\,\,\,\,\,\,\,\,\, {\rm with } \,\,\,\,\,\,\,\,\,\, 
\tilde{F}_{\ddot{n}}:=\frac{1}{2\pi} \int_{-\pi }^{-\pi }\tilde{f}(\theta ) \, e^{-i\ddot{n}\theta } \, d\theta \,\, ,
\Label{eqn:exp_fou}\end{eqnarray}
by the change of variable $x\to\theta$ and the relation  
\begin{eqnarray}
\displaystyle\tilde{F}_{\ddot{n}}=\frac{(-1)^{\ddot{n}}}{\sqrt{2\pi }}\ddot{f}_{\ddot{n}} \, \, .
\Label{eqn:equality_coeff}\end{eqnarray}

%Since the point-wise convergence of the Fourier series is guaranteed for $x\in(a,b)$ for any function which is continuous and of bounded variation in $(a,b)$, 
The correspondence introduced above provides us 
with the following theorem:
\begin{theorem}
\Label{thm:pointwise}
If a number sequence $\{f_n\}_{n=0}^\infty \in \ell ^2(\mathbb{Z}^+)$ satisfies $\, \displaystyle \lim _{N\to\infty } \Bigl\| \Bigl(\sum _{n=0}^N f_n e_n \Bigr) -f\, \Bigr\|_{(k_0)}=0$ with $f\in C^1(\mathbb{R} )$, \, then $\, \displaystyle \lim _{N\to\infty }\sum _{n=0}^N f_n e_n(x) =f(x)$ holds for any $x\in\mathbb{R}$.
\end{theorem}
{\em Proof of Theorem \ref{thm:pointwise}}
 
\rm 
From the correspondence (\ref{eqn:choice-bs}) between the unilateral orthonormal basis system $\{e_n\, |\, n\in\mathbb{Z}^+\}$ and the bilateral orthonormal basis system $\{\sqrt{\frac{1}{\pi }}\, \psi _{k_0,\ddot{n}} \, | \, \ddot{n}\in\mathbb{Z}\}$ of ${\cal H}$, it is easily shown that the coefficients $\ddot{f}_{\ddot{n}}$ $(\ddot{n}\in \mathbb{Z})$ in the expansion (\ref{eqn:exp_psi}) correspond to the coefficients $f_n$ $(n\in\mathbb{Z}^+ )$ in the expansion $\displaystyle f(x)=\sum_{n=0}^\infty f_ne_n(x)$ by the relation $f_n=\ddot{f}_{\ddot{n}_{k_0,n}}$. Hence, the condition $f\in L_{(k_0)}^2(\mathbb{R} )$ is equivalent to the conditions $\{f_n\}_{n=0}^\infty \in \ell ^2(\mathbb{Z}^+ )$, $\{\ddot{f}_{\ddot{n}}\}_{\ddot{n}=-\infty }^\infty\in \ell ^2(\mathbb{Z} )$ and $\{(-1)^{\ddot{n}}\tilde{F}_{\ddot{n}}\}_{\ddot{n}=-\infty }^\infty\in \ell ^2(\mathbb{Z} )$, where the equivalence to the last condition can be shown by (\ref{eqn:equality_coeff}). Hence, under the condition $f\in L_{(k_0)}^2(\mathbb{R} )$, the function $\tilde{f}(\theta )$ defined in (\ref{eqn:def_tildef}) for $f(x)$ in this lemma belongs to $L^2([-\pi, \pi])\subset L^1([-\pi,\pi])$. %Therefore, $\tilde{f}(\theta )$ is Lebesgue-integrable. 
This fact implies that $\displaystyle \int_{\pm\pi \mp\epsilon}^{\pm\pi} \frac{|\tilde{f}(u)|}{\theta -u}\, du <\infty $ and $\displaystyle \int_{\pm\pi \mp\epsilon}^{\pm\pi} \frac{|\tilde{f}(u)|}{\theta +u}\, du <\infty $\, if $0<\epsilon <\frac{\pi }{4}$ and $|\theta |\le \pi -2\epsilon $. 
Moreover, the function $\tilde{f}(\theta )$ is continuously differentiable once in the interval  $[-\pi+\epsilon ,\,  \pi-\epsilon ]$, % for $0 < \epsilon <\frac{\pi}{4}$, %for any $\epsilon >0$, 
from (\ref{eqn:def_tildef}) and $f\in C^1(\mathbb{R} )$. 
%Therefore, $\frac{d}{d\theta }\tilde{f}(\theta )$ is bounded in $[-\pi+\epsilon ,\,  \pi-\epsilon ]$ for any $\epsilon >0$, and hence it is of bounded variation in $(-\pi+\epsilon ,\,  \pi-\epsilon )$ for any $\epsilon >0$. 

From these facts, for $\varphi_{\theta }(t):=\tilde{f}(\theta  +t)+\tilde{f}(\theta  -t)-2\tilde{f}(\theta )$, it is easily shown that $\displaystyle \int_0^\pi \frac{|\varphi_{\theta }(t)|}{t}dt <\infty $ for any $\theta \in [-\pi+2\epsilon, \pi-2\epsilon ]$. Thus, Dini's test~\cite{enc} is passed for point-wise convergence of the Fourier series, for any $\theta \in [-\pi+2\epsilon, \pi-2\epsilon ]$.
Hence, from (\ref{eqn:psi_exp}) and (\ref{eqn:exp_fou}), 
 $\, \displaystyle \lim _{N\to\infty }\sqrt{2}\sum _{\ddot{n}=-N}^N \,(-1)^{\ddot{n}}\tilde{F}_{\ddot{n}}\,\tilde{\psi } _{k_0,\ddot{n}}(\theta ) =  \lim _{N\to\infty }\sum _{\ddot{n}=-N}^N \tilde{F}_{\ddot{n}}\,e^{i\ddot{n}\theta }=\tilde{f}(\theta )$ holds for any  $\theta\in (-\pi +2\epsilon,\,  \pi -2\epsilon )$ for any $0< \epsilon <\frac{\pi}{4}$. %because the point-wise convergence of the Fourier series is guaranteed in this interval due to the continuity and the bounded variation for $\theta \in(-\pi +\epsilon,\,  \pi -\epsilon 
%  )$. 
Moreover, 
$\, \displaystyle \lim _{N\to\infty }\sqrt{2}\sum _{\ddot{n}=-N-k_0-d}^{-N-1} (-1)^{\ddot{n}}\tilde{F}_{\ddot{n}}\,\tilde{\psi } _{k_0,\ddot{n}}(\theta ) = 0$ \,  $(d=1,2)$ for any $\theta\in (-\pi +2\epsilon,\,  \pi -2\epsilon )$, because $\{(-1)^{\ddot{n}}\tilde{F}_{\ddot{n}}\}_{\ddot{n}=-\infty }^\infty \in \ell ^2(\mathbb{Z}^+)$ (as is mentioned above) and 
\begin{eqnarray*}
\Bigl|\sum_{\ddot{n}=-N-k_0-d}^{-N-1}(-1)^{\ddot{n}}\tilde{F}_{\ddot{n}}\tilde{\psi}_{k_0,\ddot{n}}(x)\Bigr|^2 
&\le &
\Bigl( \sum_{\ddot{n}=-N-k_0-d}^{-N-1} |\tilde{\psi} _{k_0,\ddot{n}}(x)|^2 \Bigr) \, \Bigl( \sum_{\ddot{n}=-N-k_0-d}^{-N-1} |\tilde{F}_{\ddot{n}}|^2 \Bigr) \\ 
&\le & \frac{k_0+d}{(1+x^2)^{k_0+1}} \,\,  \Bigl( \sum_{\ddot{n}=-\infty }^{-N-1} |\tilde{F}_{\ddot{n}}|^2 \Bigr) . 
\end{eqnarray*}
Therefore, $ \displaystyle \lim _{N\to\infty }\sqrt{2}\sum _{\ddot{n}=-N-k_0-d}^N (-1)^{\ddot{n}}\tilde{F}_{\ddot{n}}\tilde{\psi } _{k_0,\ddot{n}}(\theta ) = \tilde{f}(\theta )$  $(d=1,2)$ for any 
$\theta\in (-\pi +2\epsilon,  \pi -2\epsilon )$. This fact and (\ref{eqn:equality_coeff}) imply that 
$\, \displaystyle \lim _{N\to\infty }\!\!\! \sum _{\ddot{n}=-N-k_0-d}^N \!\!\!  \ddot{f}_{\ddot{n}}\,\psi _{k_0,\ddot{n}}(x) =f(x)$\,  $(d=1,2)$ for any $\displaystyle x\in \left(\tan\frac{-\pi+2\epsilon }{2},  \tan\frac{\pi-2\epsilon }{2}\right)$ for any $0< \epsilon <\frac{\pi }{4}$.  Since $\displaystyle\lim_{2\epsilon \to 0+} \tan\frac{\pm\pi\mp2\epsilon }{2}=\pm\infty $,
$\, \displaystyle \lim _{N\to\infty }\sum _{\ddot{n}=-N-k_0-d}^N \ddot{f}_{\ddot{n}}\,\psi _{k_0,\ddot{n}}(x) =f(x)$\,  $(d=1,2)$ holds for any $x\in\mathbb{R} $. Since the `matching' in (\ref{eqn:choice-bs}) results in  $\displaystyle\sum _{\ddot{n}=-N-k_0-d}^N \ddot{f}_{\ddot{n}}\,\psi _{k_0,\ddot{n}} =  \sum_{n=0}^{2N+k_0+d-1} f_n e_n $ \, $(d=1,2)$ where the last equality should hold because $\{e_n\, |\, n\in\mathbb{Z}^+ \} $ is a basis system of ${\cal H}$,\, the convergences   $\, \displaystyle \lim _{N\to\infty }\sum _{n=0}^{2N} f_n e_n(x) =f(x)$ and $\, \displaystyle \lim _{N\to\infty }\sum _{n=0}^{2N+1} f_n e_n(x) =f(x)$ hold for any $x\in\mathbb{R} $, and hence  $\, \displaystyle \lim _{N\to\infty }\sum _{n=0}^N f_n e_n(x) =f(x)$ holds for any $x\in\mathbb{R} $. 
\hfill\endproof
%\end{proof_theorem_x} 

\end{document}